\documentclass[leqno, 11pt, a4paper]{amsart}

\usepackage{amsmath}
\usepackage{amssymb, latexsym, slashed, amscd}
 \usepackage[all]{xy}
 \usepackage{amsthm}
\usepackage{amsfonts}
\usepackage{mathrsfs}
\usepackage{enumerate}
\usepackage{color}
  \usepackage{comment}
\usepackage{ascmac}

 \newtheorem{definition}{Definition}[section]
 \newtheorem{theorem}[definition]{Theorem}
 \newtheorem{lemma}[definition]{Lemma}
 \newtheorem{proposition}[definition]{Proposition}
 \newtheorem{corollary}[definition]{Corollary}
  \newtheorem{conjecture}[definition]{Conjecture}

  \newtheorem{sublemma}[definition]{Sublemma}

 \newtheorem*{theorem*}{Theorem}
\newtheorem*{proposition*}{Proposition}
\newtheorem*{lemma*}{Lemma}

 \theoremstyle{remark}
 \newtheorem{example}[definition]{Example}
 \newtheorem{remark}[definition]{Remark}

\newcommand{\op}[1]{\operatorname{#1}}

\def\XXint#1#2#3{{\setbox0=\hbox{$#1{#2#3}{\int}$}
\vcenter{\hbox{$#2#3$}}\kern-.5\wd0}}

\newcommand{\C}{\ensuremath{\mathbb{C}}}

\newcommand{\R}{\ensuremath{\mathbb{R}}} 
 
\newcommand{\Z}{\ensuremath{\mathbb{Z}}}

\newcommand{\supp}{\op{supp}}

\newcommand{\coker}{\op{coker}}

\newcommand{\ba}{\begin{eqnarray}}
   \newcommand{\na}{\end{eqnarray}}

\numberwithin{equation}{section}

\begin{document}


\title{
Holomorphic curves
into infinite dimensional almost K\"ahler manifolds
and Hamiltonian dynamics}

\author{Tsuyoshi Kato}
 \address{Department of Mathematics, Kyoto University, Japan}
 \email{tkato@math.kyoto-u.ac.jp}

\begin{abstract}
We study  analysis over infinite dimensional manifolds
consisted by sequences of almost K\"ahler manifolds.
In particular 
we develop moduli theory of pseudo holomorphic curves 
into the spaces with high symmetry.

As applications, we study Hamiltonian dynamics over the
infinite dimensional manifolds, and induce some dynamical properties
of Hamiltonian diffeomorphisms on such spaces.
\end{abstract}


\maketitle

\section{Introduction}
Global  analysis of 
infinite dimensional spaces  is one of the
main subjects   in geometry.
In this paper, we  study 
 infinite dimensional  geometry and analysis from 
symplectic geometry  view point. In particular
we  introduce  a class of 
 infinite dimensional geometric spaces 
consisted by sequences of embeddings by finite dimensional manifolds.
We develop a  basic  analytic tool to perform some functional analysis
under the conditions of  high symmetry over such spaces.

 Moduli theory of pseudo holomorphic curves
is    the very powerful tool and has become one of the central theme
in symplectic geometry. From the  view point of infinite dimensional geometry,
we  study   
{\em  almost K\"ahler sequences}  $[(M_i, \omega_i, J_i)] $  which consist of
 families of   embeddings by  almost K\"ahler manifolds:
$$  (M_0 , \omega_0, J_0) \subset (M_1 , \omega_1, J_1) 
                \subset \dots \subset (M_i, \omega_i, J_i) \subset \dots $$ 
                In order to develop global analysis over these spaces,
                we introduce Sobolev spaces with Hilbert space coefficients, and then 
                we study 
                the moduli theory of pseudo holomorphic curves into such  
                infinite dimensional spaces
                 from two dimensional sphere.
  Our  construction provides with  two main ingredients.
One is Fredholm theory of the linearized maps,
where it requires two conditions;
 closedness of the range and  well definedness of their indices.
The other is non linear analysis where it also requires two conditions;
 regularity of  maps and
compactness of the moduli spaces. 
We discover that these properties are also satisfied over 
the   infinite dimensional spaces
under the    high symmetry conditions.
  Many of the   complex infinite homogeneous spaces satisfy such symmetries.

Let  $[(M_i, \omega_i, J_i)] $ be an almost K\"ahler sequence, and hence
each $(M_i, \omega_i, J_i)$ is an almost K\"ahler manifold.
A map $u: S^2 ={\bf CP}^1 \to M_i$ is called {\em holomorphic}, 
if the almost complex structures commute at each $z \in S^2$:
$$J_i \circ Tu_z = Tu_z \circ \sqrt{-1}: T_z S^2 \to T_{u(z)}M_i.$$

Let 
$E(J)_i, \ F_i \mapsto S^2 \times M_i$ be 
 vector bundles whose fibers are given by endomophisms respectively:
\begin{align*}
& E(J)_i(z,m) = \{ \phi: T_z S^2 \mapsto T_mM_i : 
   \text{ anti complex linear } \}, \\
& F_i(z,m) = \{ \phi: T_z S^2 \mapsto T_mM_i : \text{ linear} \}.
\end{align*}
Let us fix a large $l \geq 1$. Then for each $i$,  let: 
 $${\frak B}_i \equiv  {\frak B}_i (\alpha)$$ 
 be the sets of normalized and pointed
$L^2_{l+1}$ maps from $S^2$ to $M_i$, 
whose homology classes represent $\alpha \in \pi_2(M_i)$.
Then we have   two stratified Hilbert bundles over ${\frak B}_i$:
 \begin{align*}
&  {\frak E}_i = L^2_l({\frak B}_i^*(E(J)_i)) = 
\cup_{u \in {\frak B}_i}  \ \{u\} \times L^2_l(u^*(E(J)_i)), \\
&   {\frak F}_i = L^2_l({\frak B}_i^*(F_i)) 
= \cup_{u \in {\frak B}_i}   \ \{u\} \times L^2_l(u^*(F_i)).
\end{align*}
Notice that these spaces admit continuous $S^1$ actions
induced from the canonical action on ${\bf CP}^1$.

 The non linear Cauchy-Riemann operator is  given  as  sections:
$$ \bar{\partial}_i \in C^{\infty}({\frak E}_i \mapsto {\frak B}_i),  
\quad
  \bar{\partial}_i (u) = Tu + J \circ Tu \circ \sqrt{-1} .$$
$u$ is called a {\em holomorphic curve} if it satisfies the equation 
$\bar{\partial}_i(u)=0$.
The moduli space of holomorphic curves is defined by:
$$     {\frak M}(\alpha, M_i,J_i)  = 
    \{ \ u \in C^{\infty}(S^2, M_i) \cap {\frak B}_i(\alpha) 
      : \bar{\partial}_{i} (u) =0  \  \}.$$
$J$ is called {\em regular},
 if the linearized maps
  are onto at all $u \in    {\frak M}(\alpha, M_i,J_i) $ and all $i \geq 0$:
 $$D \bar{\partial}_{J}(u): T_u {\frak B}_i \mapsto ({\frak E}_i)_u$$

Take any $u \in {\frak B} \equiv \cup_{i \geq 0}  \ {\frak B}_i$, and consider  an open neighborhood
$U(u) \subset {\frak B}$. By introducing  Sobolev norms, one can make completion
$U(u) \subset \hat{U}(u)$. 
Notice that an element in $\hat{U}(u)$ cannot be realized 
by  a map into $M =\cup_{i \geq 0} \  M_i$ in general.
Let: 
 $$ \bar{\partial}_{J}: \hat{U}(u) \mapsto \hat{\frak E}|\hat{U}(u)$$
 be the extension of the CR operator on their completions.
 The differential of the operator  is not necessarily onto even if it is regular, where
 the range may not be closed.
 An almost K\"ahler sequence is said to be {\em strongly regular},
  if  the extensions are onto  at all:
 $$u \in   {\frak M}[(M_i, \omega_i, J_i)]   \equiv \cup_{i \geq 0}
  \   {\frak M}(\alpha, M_i,J_i) .$$
 This is a key property of moduli theory we develop in the infinite dimensional setting.

 \begin{theorem}\label{basic-thm}
 Let $[(M_i, \omega_i, J_i)]$ be a symmetric  K\"ahler sequence.
 
$ (1)$
If it is  regular and $\dim  \cup_{i \geq 0}    \ker D_u \bar{\partial}_i =N$ is finite,
 then it is in fact strongly regular of index $N$.

 In particular  
  $ {\frak M}[(M_i, \omega_i , J_i)] $ is a regular $N$ dimensional manifold.
  
  $(2)$  If  moreover  it is  isotropic and is regular with respect to 
  a minimal class  $\alpha \in \pi_2(M)$, then 
   the equality holds:
   $$ {\frak M}[(M_i, \omega_i , J_i)] =  {\frak M}( M_0,\omega_0, J_0) .$$
   Moreover they are compact.
    \end{theorem}
   
   As a particular case,
    ${\bf CP}^{\infty}$ with the Fubini-Study form satisfies $(1)$ and $  (2) $
     above with $N=1$.
   Moreover the moduli space is homeomorphic to $S^1$.

If $\alpha \in \pi_2(M)$ is non zero, then there is a free $S^1$ action on
$ {\frak M}[(M_i, \omega_i , J_i)]$.
In particular in the case of  ${\bf CP}^{\infty}$, the moduli space 
 is $S^1$-freely   cobordant to non zero.

\vspace{3mm}

Let us apply the above framework to Hamiltonian dynamics 
over $M= \cup_{i \geq 0} \ M_i$.
Let $f : M \to \R$ be a smooth and bounded function
which we call a bounded Hamiltonian function.
We  introduce  three  classes on such functions as
{\em pre-admissiblity}, {\em properness} and {\em connectedness}
(definition \ref{ham.fcn-def}).
The first condition  is necessary for analysis of  moduli spaces.
The second one is to induce  a Hamiltonian diffeomorphism on $M$
which is an $\R^{\infty}$ manifold rather than Hilbert one.
The third one is 
to construct an action functional invariant.

\begin{definition}
Let $\mathcal M$ be a class of diffeomorphisms on $M$, and let 
 $F: M \cong M$ be an infinitely cyclic diffeomorphism.
 
It is stably infinite cyclic in $ \mathcal M$,
 if there is a neighborhood $\mathcal{N}$ of $F$ in  $ \mathcal M$
 such that
any element $G \in \mathcal N$ is also infinitely cyclic.
\end{definition}
A simple example of an unstably infinite cyclic diffeomorphism is 
an irrational rotation on the circle.

For our study on stability of  infinite cyclicity of a Hamiltonian diffeomorphism, 
we use $S^1$ equivariant moduli theory of  holomorphic curves 
into $M = \cup_i \ M_i$.

\begin{theorem} \label{unif.cyc}
Let $[(M_i, \omega_i, J_i)]$ be an 
 almost K\"ahler sequence with 
 finite diameter of $M$.

Assume that   ${\frak M}[(M_i, \omega_i, J_i)]$
is non empty, regular, compact and $S^1$-freely cobordant to non zero
with respect to a minimal class in $\pi_2(M)$.

Let
$F: M \cong M$ be the non trivial Hamiltonian diffeomorphism induced from a bounded Hamiltonian function $f$
which is pre-admissible, proper and $2$-connected.

Then $F$ is stably  infinite cyclic in the set of Hamiltonian diffeomorphisms.
\end{theorem}

           Let us take 
$ u \in {\frak M}[(M_i,\omega_i,J_i)]$, and  $U(u) $
 be  an open  neighborhood of $u \in {\frak B} \equiv \cup_i \  {\frak B}_i$. 
 Any element in $U(u)$ is locally modelled on maps
 from ${\mathbb R}^2$ into ${\mathbb R}^{2 \infty}$
 passing through local charts.
 By use of the embedding into Hilbert space ${\mathbb R}^{2 \infty} \subset H$,
one can introduce Sobolev spaces and take completion 
 to a Hilbert manifold $\hat{U}(u)$.
Then let us put: 
$$ \hat{\frak M} [(M_i,  \omega_i, J_i)]  = 
\cup_{u \in {\frak M}[(M_i, \omega_i, J_i)] }
   \ \  \{ \ v \in  \hat{U}(u)
      : \bar{\partial}_{J} (v) =0  \  \}.$$
If $ \hat{\frak M} [(M_i,  \omega_i, J_i)]  \backslash  {\frak M} [(M_i,  \omega_i, J_i)] $ could be non empty,
then one can find some sequence $u_i \in {\frak B}_i $ 
which do not converge, but satisfy: 
$$\lim_{i \to \infty}  \  ||  \bar{\partial}_i(u_i)  || \ = \ 0.$$
Actually  we verify that such phenomenon cannot happen under some conditions.

\begin{definition}
Let  $D$ be a subset with 
${\frak M}  [(M_i,  \omega_i, J_i)] \subset D \subset {\frak B}$.

It  is 
  {\em properly compact} over  ${\frak M} [(M_i,  \omega_i, J_i)] $,
 if   any subset $\{u_k \}_k \subset  D$ with 
$\lim_{k \to \infty} ||\bar{\partial}_J(u_k)|| =0$ admits
 a subsequence which converges to some element in 
${\frak M} [(M_i,  \omega_i, J_i)] $.
\end{definition}

Let $[(M_i, \omega_i, J_i)]$ be an almost K\"ahler sequence.
Let us   consider
a bounded Hamiltonian $f: M \to [0, \infty)$ with 
 the restrictions  $f_i  =f|M_i $.
Then 
we construct a family of     functionals 
${\frak F}_i : {\frak B}_i \to {\frak E}_i$  from $f_i$,
and the (possibly non compact) cobordisms:
 \begin{align*}
 {\frak C}_i  =   \{ \ (u, \lambda) : \bar{\partial}_i(u) + \lambda {\frak F}_i(u)=0 \ \} \
  \subset  \ {\frak B}_i \times [0, \infty)
\end{align*}
with
$\partial {\frak C}_i = {\frak M}(M_i, \omega_i, J_i).$

It turns out that  there is a bound
 ${\frak C}_i \subset {\frak B}_i \times [0, \delta_0]$ for some $\delta_0 \geq 0$ 
 which is given by some topological number and is 
 independent of $i$.

Let us put:
$${\frak D}= 
\{u :   \  (u, \lambda)  \in {\frak C}_i   \text{ for some } i  \text{ and } \lambda  \}
\ \subset \ {\frak B}$$
with ${\frak C}_i(\lambda ) = \{ \ u : (u, \lambda) \in {\frak C}_i \ \}$ 
and ${\frak C}_i(0 \leq \lambda <\lambda_0)= \cup_{0 \leq \lambda <\lambda_0} \
{\frak C}(\lambda)$.

We verify the following product structure:
\begin{theorem}\label{product}
Let $[(M_i, \omega_i, J_i)]$ be an isotropic,
  quasi transitive and K\"ahler sequence.
Suppose  that the moduli space
${\frak M}[(M_i, \omega_i, J_i)]$ is non empty, finite dimensional 
and regular
with respect to a minimal class.

 If a bounded Hamiltonian $f$  is proper and pre-admissible, then  the corresponding 
${\frak D}$ 
 is properly compact.
 
 In fact there is a positive $\lambda_0 >0$ so that  the homeomorphisms:
 $$  {\frak C}_i(0 \leq \lambda <\lambda_0) \ 
\cong \ {\frak M}(M_0, \omega_0, J_0) \times [0, \lambda_0)$$
hold for all sufficiently large $i$.
\end{theorem}

In particular under the above conditions, 
 each slice  ${\frak C}_i(\mu )$ is homeomorphic to a compact space 
$ {\frak M}(M_0, \omega_0, J_0) \times \{ \mu \} $ for $0 \leq \mu < \lambda_0$.
In order to obtain such uniform product structure, we use regularity of the moduli space
and apply  the implicit function theorem.
So in a general situation,
it would be too much to expect to obtain some
 positive $\lambda_0 >0$ with the property.

On the other hand non compactness of ${\frak C}_i$ plays the essential role 
in the proof of theorem \ref{unif.cyc}.

Motivated by these phenomena, we introduce a numerical invariant of the cobordism:
 \begin{align*}
 { Cob}(f) = 
 \liminf_i  
 \sup_{\lambda \geq 0} \{ \   \lambda   :  \
  {\frak C}_i(\mu) 
 \text{ are } & \text{   non empty } \\
 & \text{ and compact  for all  } 0 \leq \mu \leq  \lambda  \}.
 \end{align*}

  \vspace{3mm}

   Let us introduce an asymptotic growth invariant of  Hamiltonian diffeomorphisms.
   Let $f: M \to [0, \infty)$ be a bounded Hamiltonian and $f_l: M_l \cong M_l$ be the 
   restrictions.  Then we obtain the Hamiltonian diffeomorphisms 
   $F_l :M_l \cong M_l$.

     Let us  introduce:
  $$\text{ As} (f)  \equiv \lim \inf_{m} \  \lim \inf_{l}
   \frac{1}{m} \   ||d(F_l)^m||_{C^0(M_l)}.$$

  These two invariants focus on different natures on functions.
 We  induce a new inequality.
 Let us fix subsets  $N_0, N_{\infty},  U \subset M$ in pre-admissibility and connectedness for $f$.
  \begin{theorem}\label{inequ}
 Let $[(M_i, \omega_i, J_i)]$ be an
  isotropic, quasi transitive and K\"ahler sequence
  with finite diameter.
 Assume  that ${\frak M}[(M_i, \omega_i, J_i)]$ is non empty, regular, 
 and $S^1$-freely cobordant to non zero
 with respect to a minimal class.

Then there is a constant $C \geq 0$  determined by $[(M_i, \omega_i, J_i)]$
and $(N_0,N_{\infty}, U)$
so that 
a uniform positivity:
$$C \  \leq \  \text{ Cob }( f )   \text{ As}( f)$$
holds for any non trivial bounded Hamiltonian function $f$
which is pre-admissible, proper and $2$-connected.
\end{theorem}
Roughly speaking this estimate saids that if $f$ is more `complicated'
so that { Cob}$(f)$ becomes smaller, 
then the norms of the derivatives of the iterations should be bigger.

\vspace{3mm}

Let us describe a possible further research direction.
Actually it is our original motivation.

The infinite projective space
${\bf CP}^{\infty}$ equipped with the homogeneous coordinate
$[z_0,z_1, \dots]$ is K\"ahler with  the Fubini-Study  metric.
We can apply moduli theory over almost K\"ahler sequences as above,
since it satisfies all analytic conditions appeared in this paper.

Our motivation of this study arose from construction of continuous
deformation of discrete  groups.
Even though ${\bf CP}^{\infty}$  is a particular space, one can immediately
see  that
its automorphism group  
Aut ${\bf CP}^{\infty}$ preserving  the K\"ahler structure
 is quite rich.
It contains any discrete group $\Gamma$ acting on a tree as a subgroup:
 $$\Gamma \subset \text{ Aut } {\bf CP}^{\infty}.$$
  The construction is quite simple. Let $V$ be the set of vertices on $T$, 
  and  assign indices arbitrarily to all the elements as
  $V = \{v_0,v_1, \dots\}$. Then the action $g \in \Gamma$ on $T$ induces
  an element in Aut ${\bf CP}^{\infty}$ by:
  $$[z_0,z_1, \dots] \to [z_{e(g(v_0))},z_{e(g(v_1))}, \dots]$$
  where $e(v_i) =i$.
Even though this element of course depends on choice of the assignment,
 it turns out that the  embedding of the group 
 $\Gamma$ above is  canonical  up to conjugacy.
 Such direction of study has been already developed in  \cite{gromov3}.

Study on
infinite groups acting on    trees 
  is  an important branch in infinite group theory, 
presented  by   Bass-Serre theory \cite{serre}  and 
 theory of  automata groups.  
 An automata group is   finitely generated
 and constructed from a Mealy automaton, which  
 act on  a rooted tree.
It has been known that  automata groups contain important  classes of 
finitely generated groups. Particular instances are given by
  intermediate growth groups which gave a counter example to
the Milnor's conjecture \cite{grigorchuk2}, construction of 
 finitely generated infinite torsion groups which solved
the Burnside problem \cite{al}, \cite{grigorchuk1}, and 
finitely generated  groups with non-uniform growth functions \cite{wilson}.

 Hamiltonian deformation of 
finitely generated groups acting on  trees
gives an
interplay between infinite group theory and symplectic geometry.
Let $F : {\bf CP}^{\infty} \cong {\bf CP}^{\infty}$ be a Hamiltonian diffeomorphism,     
and consider  deformation of an automorphism $g $ by:
$$g' = F \circ g:  {\bf CP}^{\infty} \cong {\bf CP}^{\infty}.$$
An automata group $\Gamma$  admits a canonical generating set
by their states $S=\{s_1, \dots, s_l\}$ so that
one   obtains the group  deformation of $\Gamma$:
$$\Gamma' = \text{ gen } \{ F \circ s_1, \dots, F \circ s_l \} 
\subset \text{ Diff } {\bf CP}^{\infty}.$$
Automata groups with small states have been classified 
\cite{grigorchuk-etal} (see \cite{grigorchuk-zuk}),
where many finite groups appear in the class.
One may ask whether there could exist a Hamiltonian deformation 
 which preserves finiteness, 
or more interestingly which produces  infinite torsion groups.
An immediate answer to the former is given  just by
 a    Hamiltonian rotation.

\section{Almost K\"ahler sequences}
\subsection{Stratified  local charts}
Let us introduce notations of basic open subsets.
For positive $\epsilon >0$, 
let: 
$$D^{2k}( \epsilon)  \subset {\mathbb R}^{2k}$$ be 
$\epsilon$ ball with the center $0$. We denote the $2i$ dimensional $\epsilon$ cube:
$$D_i = D^2(\epsilon) \times \dots \times D^2(\epsilon)$$
 by multiplication of $D^2(\epsilon)$ by  $i$ times.
There are canonical embeddings: 
$$D_i =D_i  \times \{0\}   \ \subset \ D_{i+1}$$ for all $i \geq 1$.
Let us put the infinite dimensional cube and disk by:
$$ D_{\infty}   \equiv  \cup_{i \geq 1} \ D_i, \quad D(\epsilon ) \equiv  \cup_{k \geq 1} \ D^{2k}(\epsilon) 
 \ \  \subset \  {\mathbb R}^{\infty}$$
 respectively. Notice   diam $D_{\infty}= \infty$.

 \subsubsection{Various norms}\label{var-norm}
 Let $H$ be the separable Hilbert space which is obtained by the 
 completion of ${\mathbb R}^{\infty}$ with the standard metric on it.
For $p= (p_0,p_1, \dots)  \in  {\mathbb R}^{\infty}$, 
let us denote by $D_{\infty}(p)  \equiv D_{\infty}+p$
and $D(\epsilon)(p) \equiv D(\epsilon) +p$  as
  the infinite dimensional cube and disk  with the center $p$ respectively.
We denote  the metric  completion by:
 $$\bar{D}_{\infty}(p) , \quad \bar{D}(\epsilon)(p) \quad   \subset  \quad H.$$
 A  neighborhood of $p \in  {\mathbb R}^{\infty}$
 is an 
open subset $p \in B \subset {\mathbb R}^{\infty}$
so that $B$ contains some ball $D(\delta)(p')$  at any $p' \in B $,
where $\delta >0$ depends on $p'$.

Let $ p \in B \subset {\mathbb R}^{\infty}$ be an open subset, and denote its closure by
 $\bar{B} \subset H$.
Let us consider a smooth and bounded function 
$f: B \mapsto {\mathbb R}$. We will regard the derivatives of $f$ at $p$
as the linear operators:
\begin{align*}
& \nabla f:  T_p {\mathbb R}^{\infty}  \ \equiv   \ \cup_{k \geq 1}  \ T_p {\mathbb R}^{2k}
 \  \mapsto \ {\mathbb R}, \\
& \nabla^2 f  : (T_p {\mathbb R}^{\infty})^{\otimes 2}  \ \mapsto \  {\mathbb R},   \\
& \dots
\end{align*}
where
$\nabla^2 (f) (v,w) = \frac{\partial^2}{\partial s \partial t} f(p+sv+tw)|_{s=t=0}$.

For $l \geq 0$, let us denote the operator norms by  $|\nabla^l f|(p)$, if
 it extends to a bounded  linear functional:
$$\nabla^l f : (T_p H)^{\otimes l}  \ \mapsto \  {\mathbb R}.$$

\begin{definition}
Let $B \subset \R^{\infty}$ be an  open subset, and $f: B \to \R$ 
be a function.

$f$ is of {\em  completely  $C^k$ bounded geometry}
at $p=(p_0, p_1, \dots) \in B$, if there is another open subset 
$p \in B' \subset B$ such that:

\vspace{3mm} 

$(1) $
 $f|B'$ extends to a continuous function 
 $f : \bar{B}' \to \R$, 
 
$(2)$
 each  differential extends continuously:
$$\nabla^l f : T_pH \mapsto {\mathbb R}$$
 for all $0 \leq l \leq k$
  (hence  $|\nabla^l f|(p) < \infty$ hold for all $l$). 
\end{definition}
 We say $f$ is of  completely  $C^k$ bounded geometry, if it is at any point
$p \in B$  and satisfies  uniformity:
$$||f||^2_{C^k(B)}  \equiv 
\sup_{p \in B} \ \sum_{0 \leq l \leq k} |\nabla^l f|^2(p) \  <  \ \infty.$$  
 $C^{\infty}$ completely bounded geometry  is just said as 
completely bounded geometry.

A pointwise operator $D$ on functions over $B$
is of completely bounded geometry at  $p \in B$, 
if $D$ extends to a smooth  operator   over $C^{\infty}(B')$ for some $p \in B' \subset B$.

It is said just as completely bounded geometry, 
if it is    at each $p \in B$ so that
 the followings are satisfied:　　

\vspace{3mm}

$(1) $ There is a constant $C$ with:
$$D: C^0(B) \mapsto C^0(B), \quad |Df|(p) \leq C | f|(p), \quad
p \in B.$$ 
We denote its  pointwise operator norm by $||D||(p)$.

\vspace{2mm}

$(2)$  For all $k$, the following norms are all finite:
$$||D||^2_{C^k(B) }  \equiv 
\sup_{p \in B} \sum_{0 \leq l \leq k}
 ||\nabla^l D||^2 (p) \ < \ \infty.$$ 
 
 $D$ is  a {\em complete isomorphism}, if it is of completely of bounded geometry.
Moreover 
  there are constants $0< c<c'$
 so that  the uniform bounds hold  for each $p \in B$:
 $$c \  \leq \ ||D||_{C^0(B)}  \ \leq \ c'.$$
 
 \begin{lemma}\label{c-iso}
 If $D$ gives  a  complete isomorphism, 
 then $D^{-1}$ is also the same.
 \end{lemma}
 
 \begin{proof}
 The identity $D \circ D^{-1} $ is completely of bounded geometry,
 since $\nabla (\text{ id }) =0$ holds.
 Then it follows from the equality
 $ 0 = \nabla(D \circ D^{-1} ) = \nabla(D) \circ D^{-1} + D \circ \nabla(D^{-1})$
 that we have the  the estimate:
 \begin{align*}
 ||\nabla( D^{-1} ) ||_{C^0(B)} & =   ||D^{-1} \circ  \nabla(D) \circ D^{-1}||_{C^0(B)} \\
 & \leq ||D^{-1}||_{C^0(B)}||  \nabla(D) ||_{C^0(B)}|| D^{-1}||_{C^0(B)}  \\
 & \leq c^{-1} c' ||  \nabla(D) ||_{C^0(B)}.
 \end{align*}
 
 We can obtain similar estimates on   higher derivatives.
  \end{proof}

For pointwisely bilinear forms, one has a parallel notion
of {\em complete nondegeneracy}.
Later we will always treat almost K\"ahler  sequences whose almost complex structures,
symplectic structures or the induced Riemannian metrics are all
 completely nondegenerate.

\begin{example} Let $D^2  \subset {\mathbb R}^2$  be
 the standard  ball  with the center $0$, and
   consider   smooth functions 
$g, h : D^2 \mapsto [0,1]$
where:
\begin{align*} 
 & g(x) = exp(-\frac{|  x|^2}{1-|  x|^2}), 
 \quad  h(x) = exp(-\frac{|  x|}{1-|  x|}).
\end{align*}
Let us prepare infinite copies of 
 $g$ and $h$, and 
let us assign indices  as $g_i, h_i: D^2_i \mapsto [0,1]$
to distinguish them from each other.
Consider functions over 
$D_{\infty} = D^2_0 \times D^2_1 \times  \dots $:
$$ G=g_0g_1g_2 \dots, \quad H= h_0h_1 h_2 \dots$$ 
by the pointwise multiplication.
Both $G$ and $H$ are  smooth on $D_{\infty}$.
Then
$G$ is of completely   bounded geometry on the unit ball with the center zero, and
 $H$ is not at any point. Actually 
$H$ is not even continuous on $\bar{D}_{\infty}$.

 For example  let us choose a point 
$p \in  \bar{D}_{\infty}$ with:
\begin{align*}
&   ||p||^2_{L^2} \equiv \sum_{i=0}^{\infty} |p_i|^2 < \infty,
   \quad
 ||p||_{L^1} \equiv \sum_{i=0}^{\infty}  |p_i| = \infty.
\end{align*}
Then clearly $H(p)=0$, but $H(0) =1$. 
\end{example}

 \begin{remark}
 Let  
 $D: C^0(B, H) \mapsto C^0(B,H)$ be  a  pointwise linear  functional,
 and assume it  gives a complete isomorphism. 
 Then its inverse  also gives a complete isomorphism.
 This is verified by  the same way as lemma \ref{c-iso}.
  \end{remark}

\subsubsection{Local charts}
Let 
$$(M_0,g_0)  \subset (M_1, g_1)  \subset \dots $$
 be embeddings of Riemannian
 manifolds with $\dim M_i = 2d_i$, and  assume
 the compatibility condition:
  $$g_{i+1}|M_i = g_i.$$
 We will denote such families by $[(M_i,g_i)]$.
 For $p,q \in M_i$, let us denote their distance in $M\equiv  \cup_{i \geq 0}  \ M_i $ by:
 $$d(p,q)  \ \equiv \ \inf_{j \geq i} \ d_j(p,q) .$$
We denote $\epsilon$ tublar neighborhood of $M_i$ by
 $U_{\epsilon}(M_i) \subset M $:
 $$U_{\epsilon}(M_i)= \{ \ m \in M :  d (m, M_i) < \epsilon \ \}.$$

Recall 
$D^{2i}(\epsilon) \subset {\mathbb R}^{2i}$ 
and $ \bar{D}(\epsilon) \subset H$
in \ref{var-norm}.
Below we regard the Riemannian metric $g = \{ g_i\}_i$ as the pointwise operator 
 over its  local charts   $T_p D(\epsilon)  = \cup_{i \geq 0} \ T_p D^{2i}(\epsilon)$ 
 for $p \in D$.
If $g$ is of completely bounded geometry,
then one can extend it to an operator  on the  Hilbert bundle
$  T \bar{D} = \sqcup_{p \in D}   T_p \bar{D}$.

Let 
$D^{i}( \epsilon ) \subset D^{i+1}( \epsilon ) \subset  
 \dots \subset  {\mathbb R}^{\infty}$ be the embeddings of $\epsilon$ disks.
\begin{definition}
A Riemannian family $\{ g_i \}_i$ is uniformly bounded,  if
 there exists  positive 
$\epsilon  >0$ such that the following properties hold.

\vspace{2mm}

$(1)$ Every point $p \in M \equiv  \cup_{i \geq 0} \ M_i$ admits
a  stratified local chart:
\begin{align*}
&  
 \varphi(p) :D(\epsilon) \equiv \cup_i \ D^{i}(\epsilon)
     \hookrightarrow  M, \\
& \varphi(p)_i \equiv \varphi(p)|D^{e_i}(\epsilon)  \mapsto M_i
\end{align*}
with $\varphi(p)(0) =p$ and
 $e_i = \dim M_i$.

 
\vspace{2mm}

     $(2)$ With respect to $\varphi(p)$,
the induced  Riemannian metric $g = \{ g_i\}_i$ is 
 of completely bounded geometry so that 
 for each $l >0$, there is a constant $C(l) \geq 0$
   independent of $p$
such that  the estimate holds:
$$\sup_{p \in M} \ \sup_{m \in   D(\epsilon)} \
\sum_{0 \leq k \leq l}
|\nabla^k ( \varphi(p)^*g) |(m) \ \leq \ C(l) \qquad (*)$$

$(3)$ There is an increasing and proper function
$h:(0, \infty) \to (0, \infty)$
so that
the uniformly bounded distance property:
$$d(p,q) \geq h(d_i(p,q))$$
holds for any $i$ and  $p,q \in M_i$, 
 where $d_i$ and $d$ are the distances  on $ M_i$ and on $M $ 
respectively.
\end{definition}

We say that the stratified local chart as above  is a
{\em complete local chart}.
Also the above family $\{ (p, \varphi(p))\}$ 
is called a {\em uniformly bounded covering}.
Later on uniform implies independence of choice of points as above.

\vspace{2mm}

Let $f: M=  \cup_{i \geq 0} \ M_i \mapsto {\mathbb R}$ be a bounded  function and:
 $$ \varphi(p)^*(f) : D(\epsilon) \to {\mathbb R}$$
  be a family of  the induced functions with respect to a uniformly bounded covering.
We say that
$f$ is of {\em completely $C^k$-bounded geometry}, if 
they satisfy the estimate:
$$||f||_{C^k(M)} \equiv
\sup_{p \in M} ||\varphi(p)^*(f)||_{C^k(D(\epsilon))}  \ \leq \ C_k$$
for  some constants $C_k$ which are 
  independent of $p \in M$.

It is just of completely bounded geometry, if it is completely 
of $C^k$-bounded geometry
for
 all $k =0,1,2 \dots$


\begin{lemma}\label{Kl}
Let $[(M_i,g_i)]$ be  a uniformly bounded
Riemannian family  with
$\{  (p, \varphi(p))  \}$ and  $\epsilon >0$ as above.
Then the exponential map:
 $$\exp_p: \bar{D}(\epsilon') \mapsto \bar{D}(\epsilon)$$
  exists and is smooth
  for some $ \epsilon' >0  $,
 with respect to   the induced Riamannian metrics $\varphi(p)^*(g)$.
\end{lemma}

\begin{proof}
For a proof, see [Kl] ($p 57$, $p72$).
Notice that the geodesic coordinate does not preserve the stratifications
in general.
\end{proof}

   Let
       $f_n, g : M  \to { \mathbb R}$ be a  family of bounded  functions
       for  $n=0,1,2,  \dots$.
      We say that $\{f_n\}_n$ {\em converges weakly} to $g$ in $C^l$,
      if  the restrictions:
      $$f_n|M_k \to g|M_k$$
      converge in $C^l$
     for all $k=0,1,2, \dots$.

   \begin{lemma}
   Let $[(M_i, g_i)]$ be a uniformly bounded Riemannian family
   such that each $M_k$ is compact.
   Let
       $f_n : M  \to { \mathbb R}$ be a family of  functions 
       such that  $C^{l+1}$ norms are 
       uniformly bounded:
       $$||f_n||_{C^{l+1}(M) } \leq C(l+1)$$
         for  $n=0,1,2,  \dots$
       
    Then a subsequence $f_{n_j}$ weakly converges in $C^l$ to a function
    $g:  M \to { \mathbb R}$ of completely $C^l$-bounded geometry.
       \end{lemma}
       
       \begin{proof}
By the condition, the restrictions $\{f_n|M_k\}_n$ satisfy uniformity
of $C^{l+1}$ norms $||f_n||_{C^{l+1}(M_k)} \leq C(l+1)$.

It follows from  Rellich lemma that 
there is some $C^l$ function $g_1: M_1 \to {\mathbb R}$ 
so that  a subsequence $\{f_{n(i)}|M_1\}_i$ converges
to $g_1$ in $C^l(M_1)$.

By the same way 
there is some $C^l$ function $g_2: M_2 \to {\mathbb R}$ 
so that  a subsequence $\{f_{n(i,2)}|M_2\}_i$ converges
to $g_2$ in $C^l(M_2)$ for 
another subsequence $\{n(i,2) \}_i \subset \{n(i)\}_i$.
Clearly $g_2|M_1=g_1$ holds.

By choosing subsequences successively, 
$\{f_{n(i,k)}|M_k\}_i$ converge to some  $C^l$ function $g_k: M_k \to {\mathbb R}$
with $g_k|M_{k-1}=g_{k-1}$. 
These satisfy uniformity of $C^l$ norms $||g_k||_{C^l(M_k)} \leq c < \infty$.

Let  $g : M \to {\mathbb R}$ be a  function
defined by $g|M_k \equiv g_k$.
Then the subsequence $\{f_{n(i,i)}\}_i$ converges weakly to $g$ in $C^l$.
\end{proof}

\subsection{Almost K\"ahler sequence}
Let $(M, \omega , J)$ be a finite dimensional symplectic manifold
 equipped with a compatible almost
complex structure so that: 
 $$g( \quad , \quad ) = \omega( \quad , J \quad )$$ gives
a Riemannian metric on $M$.
Such a manifold is called  an {\em almost K\"ahler manifold}.

Let $ (M_0,\omega_0,J_0) \subset (M_1, \omega_1, J_1) 
      \subset \dots \subset  (M_i, \omega_i, J_i) \subset \dots$
be  infinite embeddings of almost K\"ahler manifolds.
  If one denotes the inclusion by 
  $I(i): M_i \hookrightarrow M_{i+1}$,
 then it implies that 
$\{ I(i)\}_i $ gives a family of holomorphic embeddings:
$$J_{i+1} \circ I(i)_* = I(i)_* \circ J_i$$
 and
 the symplectic forms are  given by 
the restrictions as:
   $$I(i)^*(\omega_{i+1})=\omega_i.$$

 Suppose  $\dim M_k=2d_k$, and let
    $U_{\epsilon}(M_i) \subset M \equiv  \cup_j M_j $
be   $\epsilon$ tubular neighborhoods of $M_i$. Let:
 $$\tilde{\pi}_k : D(\epsilon) = \cup_{j \geq 1} \ D^{2j}(\epsilon) \to D^{2d_k}(\epsilon)$$
be the standard projections.

\begin{definition}\label{almKseq}
An almost   K\"ahler sequence $[(M_i,\omega_i,J_i)]$
  consists of a family of   embeddings  by  almost K\"ahler manifolds:
$$ (M_0,\omega_0,J_0) \subset (M_1, \omega_1, J_1) 
      \subset \dots \subset  (M_i, \omega_i, J_i) \subset \dots$$
   and a positive  $\epsilon >0$
    so that    it admits 
   $\epsilon$
    uniformly bounded coverings $\{(p, \varphi(p))\}$ at all $p \in M$,
    which satisfy the followings:

\vspace{2mm}

$(1)$  $\varphi(p)^*(\omega)$ and 
  $\varphi(p)^* (J)$ are both uniformly 
 of completely nondegenerate on  $D(\epsilon)$.

\vspace{2mm}

$(2)$ The induced symplectic form is standard at $p$:
$$\varphi(p)^*(\omega)|_p  = 
    \frac{\sqrt{-1}}{2}
    \sum_{i=0}^{\infty}  \   dw_i \wedge d\bar{w}_i $$
where $(w_1, \dots, w_i)$ are the coordinates
 on $D^{2i}(\epsilon) \subset {\C}^i$.
 
 \vspace{2mm}

$ (3)$
There are     families of holomorphic  maps:
$$\pi_k : U_{\epsilon}(M_k) \mapsto M_k$$
such that 
the compatibility condition:
$$ \pi_k |  M_k  = \text{ id }, \quad 
  \pi_k(\varphi(p)(x) ) = \varphi(p)(\tilde{\pi}_k(x)) $$
  holds at
  any $p \in M_k$ and 
   any 
 $
x \in D(\epsilon)$.
\end{definition}
A complete local chart with  the properties  $(1) (2) (3)$  above, 
is called a  {\em complete almost K\"ahler chart}.

An almost K\"ahler data $\{(\omega_i, J_i)\}$ gives a uniformly bounded and
 compatible family of  Riemannian metrics on $\{ M_i\}_i$.
Notice that the equalities $< v,u> = <(\pi_k)_*(v),u>$ hold
 for $u \in T_p M_k$ and $v \in T_p U_{\epsilon}(M_k)$
 with respect to the induced Riemannian metric.

\vspace{2mm}

        Later on, we fix   a uniformly bounded covering
by  complete almost K\"ahler charts.
    
    \vspace{2mm}

We say that $[(M_i, \omega_i, J_i)]$ is  a {\em K\"ahler sequence},
if it is an  almost K\"ahler sequence  consisted by
a  uniformly bounded covering 
by  holomorphic
complete  K\"ahler charts  $ \varphi(p) $ at  all points $p$,
where we equip with  the standard complex structure on $D(\epsilon)$
(see [GH] $p107$).

Let  $f: M=\cup_{i \geq 0} \ M_i \to {\mathbb R}$ be a bounded function
on  an almost K\"ahler sequence.
 We say  that $f$  is a {\em bounded Hamiltonian function},
 if it is of completely bounded geometry.

\vspace{3mm}

Let $(M, \omega)$ be a finite dimensional symplectic manifold. 
The following  facts are  well known ([G1]):

\vspace{3mm}

(1) there exist compatible almost complex structures, and 

\vspace{2mm}

(2) the space of compatible almost complex structures
is contractible.
\vspace{3mm} \\
In our infinite dimensional situation, the condition (1) depends on the spaces,
but the same  thing holds for (2).

\begin{lemma}
Let $[(M_i, \omega_i)]$ be   a  symplectic sequence.
Suppose there exists a family of  compatible almost complex structures
$\{ J_i\}_i$ so that $[(M_i, \omega_i, J_i)]$ consists of an almost K\"ahler sequence
with respect to a uniformly bounded covering $\{(p, \varphi(p))\}$.
Then the space of such family:
\begin{align*}
{\frak J}([(M_i, \omega_i)]) =&  \{  \ \{ J'_i \}_i :  [(M_i, \omega_i, J'_i)] : \\
& \text{ almost K\"ahler  sequence with respect to } \{(p, \varphi(p))\}  \   \}
\end{align*}
is contractible.
\end{lemma}

\begin{proof}
We follow a well known argument in the  finite dimensional case.

Let us choose a reference family of almost complex structures $\{ J_i^0 \}_i$.
Take another one $\{ J_i^1\}_i$. Let us connect these by a compatible family of almost complex structures
$\{ J_i^t\}_i$, $t \in [0,1]$. For $\alpha = 0$ or $1$,
 let us put $h_i^{\alpha}( \quad , \quad ) = \omega_i(\quad , J_i^{\alpha} \quad)$.
Then $\{h_i^{\alpha}  \}_i$ gives a family of Riemannian metrics. 
Moreover each $J^{\alpha}_i$
is uniquely determined by $h_i^{\alpha}$. 
Let us consider a smooth family of Riemannian metrics:
$$h_i^t = h_i^0 + t( h_i^1 - h_i^0).$$
For each $i$, there exists a unique and smooth family of compatible almost complex structures
$J_i^t$, $t \in [0,1]$ satisfying $h_i^t (\quad , \quad) = \omega_i(\quad,  J_i^t \quad)$.
 
Let us choose 
 a complete almost K\"ahler  chart at  $p \in M_i \subset M_{i+1}$:
 $$\omega_i = \sum_{j \leq i}  \ dp_j \wedge dq_j, 
        \quad \omega_{i+1} = \sum_{j \leq i+1}  \ dp_j \wedge dq_j
\quad  \text{ at } p$$
  and denote the  local projections by $\pi'_i: D^{2d_{i+1}}(\epsilon) \mapsto D^{2d_i}(\epsilon)$.
  Let us   check  the compatibility condition
$J_{i+1}^t \circ \pi_i' = \pi'_i \circ J_i^t$  at $p$ and for each $t$.
Let us take  $v_i \in T_p M_i$. Then:
\begin{align*} 
 \omega_{i+1}(\quad ,   J_{i+1}^t v_i)  
       &    = h_{i+1}^0(\quad , v_i) + t( h_{i+1}^1 (\quad , v_i) - h_{i+1}^0(\quad , v_i)) \\
&     = \omega_{i+1}(\quad , J_{i+1}^0  v_i) + t\{ \omega_{i+1}(\quad , J_{i+1}^1  v_i) -
                                \omega_{i+1}(\quad , J_{i+1}^0  v_i) \} \\
 & = \omega_{i+1}(\quad, J_i^0v_i) + t \{    \omega_{i+1}(\quad, J_i^1 v_i) -
                                                            \omega_{i+1}(\quad, J_i^0 v_i)   \} \\
                                                            & =  \omega_{i+1}(\quad ,   J_i^t v_i)  \\
  & = \omega_i(\pi'_i  \quad , J_i^0 v_i) + t \{    \omega_i( \pi'_i \quad, J_i^1 v_i) -
                                                            \omega_i( \pi'_i \quad, J_i^0 v_i)   \} \\
  & = \omega_i( \pi'_i \quad , J_i^t v_i).
\end{align*}
The fourth equality implies the 
the  compatibility condition.

Moreover  the following equalities hold
 from the equality between the first and the last above:
\begin{align*}
&    \omega_i( \pi'_i  J_{i+1}^t(w), J_i^t \pi'_i(v)) =
\omega_{i+1}(J_{i+1}^t(w) ,   J_{i+1}^t \pi'_i(v)) \\
& = \omega_{i+1}(w ,    \pi'_i(v)) =  \omega_i(\pi_i'(w) ,    \pi'_i(v))
=  \omega_i(J_i^t(\pi_i'(w)) ,    J_i^t(\pi'_i(v))).
\end{align*}
This implies the equality: 
$$ \pi'_i  \circ J_{i+1}^t =    J_i^t \circ \pi_i'$$
and so $\pi_i'$ is holomorphic with respect to $J^t$.
\end{proof}

\subsubsection{Embeddings of almost K\"ahler sequences}
Let  $[(M_i, \omega_i, J_i) ]$ 
be an almost K\"ahler sequence equipped with 
 complete local charts $\varphi(p) :  D(\epsilon) = \cup_{s \geq 1} D^{2s}(\epsilon)
\hookrightarrow \text{ im } \varphi(p) \subset  M$ for all $p \in M$.

Let us say that $[(M_i', \omega_i', J_i') ]$ is 
 formally embeddable  into   $[(M_i, \omega_i, J_i) ]$, 
if there are 
  subindices $\{ k(i)\}_i$  and 
  compatible embeddings between almost K\"ahler manifolds:
 $$ I_i: (M_i' , \omega_i', J_i') \hookrightarrow (M_{k(i)}, \omega_{k(i)}, J_{k(i)}).$$

  \begin{example}
  $(1)$
Let us fix $p \geq 1$ and consider the canonical embeddings of the Grassmannians
$Gr_{p,q} \hookrightarrow Gr_{p,q+1}$
which embed  each $p$ plane $L \subset {\C}^{p+q} \subset {\C}^{p+q+1}$.
These admit the canonical and compatible K\"ahler forms, and 
the direct limits $Gr_p \equiv \lim_{q \to \infty} Gr_{p,q}$ consiste of the K\"ahler sequences.

Let us consider the Pl\"ucker embedding
$Gr_{p,q} \hookrightarrow {\bf CP}^N$,
where $N=N(p,q)=
\begin{pmatrix}
&p+q \\
& p
\end{pmatrix} -1$, which associate each $p$ plane $L \subset {\C}^{p+q}$
and its basis $\{v_1, \dots, v_p\} $ to the complex line
$[v_1 \wedge \dots \wedge v_p] \in {\bf CP}^N$.

It is well known that these embeddings preserve the canonical K\"ahler forms, and so 
they give the formal embedding of the K\"ahler sequences:
$$I: [Gr_{p,q}] \hookrightarrow [{\bf CP}^n]$$
where $(M_i,\omega_i, J_i) = Gr_{p,i}$ and 
$(M_i',\omega_i', J_i')={\bf CP}^i$ with $k(i) =N(p,i)$.

Moreover the Schubert calculus verifies the isomorphisms:
$$I_* : H_2(Gr_{p,q}; {\Z}) \cong H_2({\bf CP}^N; {\Z}) \cong {\Z}.$$

\vspace{3mm} 

$(2)$
Let us consider the {\em Veronese maps} defined as follows.
Let us introduce the lexicographic order for
 two indices $(i_0, \dots, i_n)$ and $(j_0, \dots, j_p)$.

Let us fix $m \in\{1,2, \dots\}$, and 
take ${\bf CP}^n$ with the homogeneous coordinate $[z_0, \dots,z_n]$.
For $N = \begin{pmatrix}
& n+m \\
& m
\end{pmatrix}
-1$, we define the Veronese map:
\begin{align*}
& v_m : {\bf CP}^n \hookrightarrow {\bf CP}^N, \\
& v_m([z_0, \dots,z_n]) = \{ z_0^{i_0}, \dots z_n^{i_n} : \sum_{l=0}^n i_l =m\} .
\end{align*}

With $n_1=1$, let us  define  numbers inductively by
$n_{i+1} = 
\begin{pmatrix}
& n_i+m \\
& m
\end{pmatrix}
-1$.

Now we have two different embeddings:
$${\bf CP}^{n_i} \subset_v {\bf CP}^{n_{i+1}}, 
\quad {\bf CP}^{n_i} \subset {\bf CP}^{n_{i+1}}$$
where the first is the given by the Veronese map and the second is by the canonical embedding.

\begin{lemma} The following diagram commutes:
$$
\begin{matrix}
{\bf CP}^{n_i} & \subset_v & {\bf CP}^{n_{i+1}} \\
\cap && \cap \\
{\bf CP}^{n_{i+1}} & \subset_v & {\bf CP}^{n_{i+2}}
\end{matrix}$$
\end{lemma}

\begin{proof}
This follows since we have used  the lexicographic order
for the coordinates. 
\end{proof}

\vspace{3mm}

\begin{corollary}
There is a canonical embeddings of ${\bf CP}^{\infty}$ to itself:
$$v_m : {\bf CP}^{\infty} \subset_v {\bf CP}^{\infty}$$
of degree $m$, so that the restrictions are given by the Veronese maps.
\end{corollary} 

\begin{remark}
We have the Veronese  sequence by the
embeddings by the iterations of the Veronese maps:
$${\bf CP}^{n_1} \subset_v {\bf CP}^{n_2} \subset_v  \dots \subset_v {\bf CP}^{n_l} \subset 
\dots \subset  V \equiv  \cup_i {\bf CP}^{n_i}  .$$
This is not K\"ahler sequence, since the degree grows unboundedly in the total space.
Study  of this embeddings will require much harder analysis.
\end{remark}
\end{example}

\subsection{Symmetric almost K\"ahler sequence}\label{sym}
Let us  introduce geometric classes of almost K\"ahler sequences.
    Their symmetric properties  allow us to analyze global 
    structure of  holomorphic maps. 

Recall the family of holomorphic maps
$\pi_k : U_{\epsilon}(M_k) \mapsto M_k$ in definition  \ref{almKseq}.
   
\begin{definition}
An almost K\"ahler sequence
$[(M_i, \omega_i, J_i)]$
 is   {\em symmetric}, if for each  $k \geq 0$
 the followings hold:

$(1)$ For each  $i \geq k+1$,
there are
  families of almost K\"ahler submanifolds:  
  $$M_k \subset W_i \subset M_i$$ 
  with $W_{k+1} = M_{k+1}$,
    and
  isomorphisms which preserve $M_k$:
$$P_i: \{(M, M_k), \omega, J\}  \cong \{(M,M_k), \omega, J\}.$$

$(2)$   $P_i$  transform $W_i$ to $M_{k+1}$ as:
$$P_i: (W_i, \omega_i|W_i, J_i|W_i) \cong (M_{k+1},\omega_{k+1}, J_{k+1})$$
such that  at any $p \in M_k$:
$$D_k :   TM_k \oplus_{i \geq l}  N_{i,k}  \cong T M |M_k $$
 gives a complete isomorphism over $M_k$ (see lemma \ref{c-iso}),
where:
\begin{align*}
& N_{i,k}= (P_i^{-1})_* [  \ (\text{ Ker } ( \pi_k)_* \cap TM_{k+1}) |M_k) \ ] \\
& D_k = \text{ id } \oplus (P_{k+1})_* \oplus (P_{k+2})_* \oplus  \dots 
\end{align*}
\end{definition}
If all these properties hold by use of complex structure,
then we say that it is a symmetric K\"ahler sequence.

\vspace{3mm}

Suppose  $[(M_i,\omega_i, J_i)]$ is a symmetric K\"ahler sequence.
It is {\em isotropic}, 
if there are families of parametrized  isomorphisms  for each  $0 \leq t \leq 1$:
$$P^t_i: \{(M, M_k), \omega, J\}  \cong \{(M,M_k), \omega, J\}$$
 with:
$$P_i^0\equiv  \text{ id }, \quad P^1_i =P_i.$$

\begin{example}
$(1)$ Let $(X,\omega, J)$ and $(Y, \tau, I)$ be two almost K\"ahler manifolds, 
and choose a base point $y_0 \in Y$. 
Let us consider the products:
$$(X \times Y_1 \times Y_2 \times \dots, \ \omega+\tau_1 +\tau_2 + \dots, \ J \oplus I_1 \oplus I_2 \oplus  \dots)$$
where all $(Y_i,\tau_i,I_i)$ are the same $(Y,\tau,I)$, and 
 we embed $X \times Y_1 \subset X \times Y_1 \times Y_2$ 
by identifying $ X \times  Y = X \times Y \times \{y_0\}$ and similar for others.

The  infinite product sequence admits symmetric structure  by
 choosing:
  $$M_k = X \times Y_1 \times \dots \times Y_k , \quad 
  W_i =M_k  \times \{ y_0 \} \times  \dots \times \{y_0\} \times Y_i.$$
 $P_i$ are given by the obvious exchange of the coordinates.
 
 \vspace{2mm}

$(2)$ Let $ M$ be a complex manifold,
and take any   holomorphic curve $u: {\bf CP}^1 \mapsto M$.
Then the holomorphic vector bundle $u^*(TM) \mapsto {\bf CP}^1$
splits as the direct sum of  holomorphic line bundles.
This fact can be regarded as
 `infinitesimal  symmetric  property' (see [OSS]).
  
\vspace{2mm}

$(3)$
The infinite complex projective space: 
$$[({\bf CP}^i, \omega_i)] = {\bf CP}^1 
\subset {\bf CP}^2 \subset \dots \subset {\bf CP}^n \subset \dots
\subset  {\bf CP}^{\infty}$$
with the Fubini Study form is an isotropic  symmetric K\"ahler sequence, and we 
 denote it  by ${\bf CP}^{\infty} \equiv \cup_{i \geq 1} \ {\bf CP}^i$.
There are standard charts ${\C}^i \subset {\bf CP}^i$ and $\omega_i$ can be expressed as:
$$\omega_i |{\C}^i =  \frac{ \sqrt{-1}}{2}
    \Bigl[  \ \frac{ \sum_l dw_l \wedge d\bar{w}_l }{(1+ w\bar{w})}
     - \frac{ ( \sum_l \bar{w}_ldw_l ) \wedge 
 (\sum_l w_l d\bar{w}_l)}{(1+ w\bar{w})^{2}} \  \Bigr]$$
where $w=(w_1, \dots,w_i)$ are the coordinates on ${\C}^i$.
The family $\{ \omega_i \equiv \omega|D^{2i}\}_i$
is completely non degenerate, where $D^{2i} \subset {\C}^i$
are the unit balls.
In order to obtain another charts at  any $p \in {\bf CP}^i$,
one can use any constant unitary matrix $U \in Mat_{i+1}({\C})$
with $U([1,0, \dots, 0]) = p \in {\bf CP}^i $.

Let $U_{\epsilon}({\bf CP}^k) \subset   {\bf CP}^{\infty}$
be $\epsilon$ tublar neighborhood. Then there are natural projections
$\pi_k : U_{\epsilon}({\bf CP}^k) \mapsto {\bf CP}^k$:
$$ \pi_k([z_0, \dots, z_k,z_{k+1},  \dots]) 
= [z_0, \dots, z_k,0 , \dots].  $$

Let us put
$M_k ={\bf CP}^k$ and $W_i$ by:
$$W_i =\{[z_0: \dots: z_k: 0  \dots : 0 : z_i: 0 : 0 : \dots ] \in {\bf CP}^{\infty} \}$$
for all $i \geq  k+1$ with  $M_{k+1}= {\bf CP}^{k+1}$.
$P_i :  W_i \cong {\bf CP}^{k+1}$ are given just by  exchange of  the coordinates:
$$[z_0: \dots: z_k: 0: \dots : 0:  z_i: 0 \dots] \to [z_0: \dots: z_k : z_i: 0: \dots].$$
This is isotropic, by putting:
\begin{align*}
 P_i^t( [z_0:  \dots: & z_k:  \dots]) 
=  [z_0: \dots:z_k  
 : \cos  \frac{\pi t}{2}  z_{k+1}+\sin \frac{\pi t }{2} z_i : \\ 
 & z_{k+2}:  \dots : z_{i-1}: 
 - \sin \frac{\pi t}{2}  z_{k+1}+\cos \frac{\pi t}{2} z_i : z_{i+1}: \dots].
\end{align*}

(4)
There are several variants. For example one can change ${\C}$ by ${\mathbb H}$.
For others, let us consider the Grassmannians:
$$Gr_{r,n}({\C}) = \{ H \subset {\C}^{r+n} \ ; \ H : r \text{ dimensional } {\C} \text{ vector subspaces } \}.$$
One can canonically embed as $H \subset {\C}^{r+n+1}$, and by taking the direct limit,
one obtains the K\"ahler sequence
$Gr_r({\C})= \lim_{n \to \infty} Gr_{r,n}({\C})$
equipped with the standard K\"ahler structure.

This space also admits isotropic and symmetric structure. Let us put:
$${\C}^{k,i}=\{ (z_1, \dots, z_k, 0, \dots, 0, z_{k+i}) : z_j  \in {\C} \} \subset  {\C}^{k+i}$$
and choose
$M_k =Gr_{r,k}$ and $W_i\equiv W^r_{k,i}$ are consisted by all elements of the form:
$$W_{k,i}^r =\{H  \subset {\C}^{k+r,i}   \ ; \ H : r \text{ dimensional } {\C} \text{ vector subspaces } \}.$$
 The required isomorphisms and isotropies
can be obtained by the same way as $(3)$.
\end{example}

\begin{lemma}\label{decomp}
Let $[(M_i, \omega_i, J_i)]$ be a symmetric   almost K\"ahler sequence.

Then there is a  bundle $N \to M_k$ 
so that 
a uniformly complete   isomorphism:
$$TM |M_k \cong TM_k \oplus (N \otimes {\mathbb R}^{\infty})$$
exists with respect to a uniformly bounded covering on $M$.
\end{lemma}

\begin{proof}
Let us put:
$$N =( \text{ Ker } (\pi_k)_* \cap TM_{k+1}  )|M_k.$$
There is a  holomorphic isomorphism $TM_{k+1}|M_k \cong TM_k \oplus N$.
Then the conclusion follows by use of the  family of 
 isomorphisms of  the tangent bundles  for al $i \geq k+1$:
$$ TM_k \oplus N \cong TW_i|M_k , \quad 
 (v,w) \to (v, (P_i^{-1})_*(w)).$$
\end{proof}

\subsubsection{Quasi transitivity}
Let $[(M_i, \omega_i, J_i)]$ be an almost K\"ahler sequence.
We  say  $[(M_i, \omega_i, J_i)]$ is {\em quasi transitive}, if
 for any  $N  >0$,
 there is  $k = k(N) $ such that  for any  points 
$p_0, \dots, p_{N-1} \in M \equiv \cup_{i \geq 0} \ M_i$,
 there is an automorphism of the almost K\"ahler
sequence $A: ((M, M_0), \omega , J) \cong ((M,M_0), \omega, J)$
which preserves $M_0$ and:
$$A(p_i) \in M_k$$
hold  for 
 all $0 \leq i \leq N-1$.

\begin{lemma}
The infinite  projective space
$[({\bf CP}^i, \omega_i, J_i)]$ is
quasi transitive. 
\end{lemma}

\begin{proof}
Let us construct automorphisms $A^i: ({\bf CP}^{\infty},{\bf CP}^{l_i}) \cong ({\bf CP}^{\infty},{\bf CP}^{l_i})$
inductively so that they satisfy the followings:
$$  A^i(p_i) \in {\bf CP}^{l_i}, \quad A^i|{\bf CP}^{l_j} = id$$
for all $j < i$.

Let us embed ${\bf CP}^{\infty} \hookrightarrow {\bf CP}^{\infty}$
by $[z_0,z_1, \dots] \to [0,0,z_0,z_1,\dots]$.
Then  diag $(1,1, A)$ is the required automorphism which preserves ${\bf CP}^1$,
where
$A \equiv A^{N-1} \circ A^{N-2} \circ \dots \circ A^0$
with $k  =l_{N-1}$.

Let us choose any $p_0=[z_0,z_1, \dots] \in  {\bf CP}^L \subset {\bf CP}^{\infty}$.
Firstly let us move $p_0$ to $ [1,0,0, \dots]$ by
a unitary  automorphism $A^0 \in U(L+1) \subset $ Aut ${\bf CP}^{\infty}$.

Let us consider $u_1= A^0(p_1) \in  {\bf C}P^{\infty}$.
 We put $A^1=$ id, if $u_1 \in {\bf CP}^1$.
Suppose $u_1 =[u_1^0,u_1^1, \dots] \not\in {\bf CP}^1$.
Then   $(u_1^1 , u_1^2, \dots)$ is non zero and so defines an element  
 in $ {\bf CP}^{\infty}$.
Let us choose another unitary automorphism $T_1$ with
$T_1([u_1^1,u_1^2, \dots]) =[1, 0, \dots]$. 
Then we put $A^1 =$ diag $(1, T_1)$.

Let us consider $u_2= A^1 \circ A^0(p_2) \in  {\bf C}P^{\infty}$.
 We put $A^2=$ id, if $u_2 \in {\bf CP}^2$.
Suppose $u_2 =[u_2^0,u_2^1, \dots] \not\in {\bf CP}^2$.
Then   $(u_2^2 , u_2^3, \dots)$  defines an element   in $ {\bf CP}^{\infty}$.
By  another unitary automorphism $T_2$ with
$T_2([u_2^2,u_2^3, \dots]) =[1, 0, \dots]$ 
Then we put $A^2 =$ diag $(1,1,  T_2)$.

By the same way one can inductively construct $A^3 , \dots, A^{N-1}$.
\end{proof}

 A similar argument can be used to verify that 
the infinite Grassmannians $Gr_N({\C}) =\lim_{L \to \infty} Gr_{N,L}$
also satisfy  quasi transitivity.

\subsubsection{Minimality}
  Let $[(M_i, \omega_i,  J_i)]$ be  an almost K\"ahler sequence. 
Let us introduce its  invariant  ([HV]):
 \begin{align*}
 m([(M_i, & \omega_i, J_i)]) 
       = \inf  \{ \ <\omega, u> ; \\
       & u: S^2 \mapsto M \equiv \cup_{i \geq 0} \ M_i 
                                                   : \text{ non constant holomorphic  curves} \ \}.
 \end{align*}
 By restriction to  the symplectic  sequence,
 one obtains another  invariant:
$$m([(M_i, \omega_i)]) = \text{ inf}_{\alpha} 
         \{  \ <\omega, \alpha> \ >0 : \alpha: S^2 \mapsto \cup_i M_i \ \}.$$
 We  say      $[(M_i, \omega_i,  J_i)]$  is {\em minimal}, if both  the equality 
 and positivity hold:
$$m([(M_i, \omega_i)]) = m([(M_i, \omega_i , J_i)]) \ > \ 0 .$$

Later on we assume that a minimal class can be represented as a map
$\alpha: S^2 \to M_0$. Actually in our arguments later, we can just shift indices
of stratification of manifolds so that this condition is satisfied.

\begin{example}
$(1)$
 Notice that if  $[(M_i, \omega_i, J_i)]$ satisfies  
 $\pi_2(\cup_i M_i) / $ Tor $\cong \Z$  of  rank $1$,
 then minimality is equivalent to existence of
non constant holomorphic curves representing a generator of $ \in  \pi_2 / $ Tor.

  The Fubini Study form on ${\bf CP}^n$ with the standard complex structure
 gives $\pi_2$ rank one minimal data $(\omega, J)$ with $m =  \pi$.

 $(2)$ Let $({\bf CP}^1, \omega, J)$ be the standard curve and
 $[(M_i, \omega_i, J_i)]$ be minimal. Then the product
 $[((M_i \times {\bf CP}^1, \omega_i+ \omega, J_i \oplus J)]$
 is also minimal.

$(3)$ Suppose $[(M_i, \omega_i, J_i)]$ is algebraic with 
each $\omega_i \in H^2(M_i ; \Z)$.
Then it is  minimal, if any generating elements  in $H_2(M; \Z)$
can be represented by  some holomorphic curves.
In particular it is the case when it  is simply connected, algebraic,
and any generating elements  in $\pi_2(M)$
can be represented by  some holomorphic curves.
\end{example}

\section{ Moduli spaces  of   holomorphic curves}
We study theory of holomorphic curves into almost K\"ahler sequences.
In particular we develop  analytic tools to 
construct finite dimensional moduli spaces over   sequences which satisfy 
some symmetric properties.

\subsection{Finite dimensional preliminaries}\label{f.d.pre}
 We recall basic materials on moduli theory of holomorphic curves into finite dimensional 
  symplectic manifolds.
  Most  of the contents   have appeared in  \cite{hofer and viterbo}.
  Based on the finite dimensional setting, 
 we  formulate Sobolev spaces over the infinite dimensional spaces  
 $M= \cup_{i \geq 0} \ M_i$.

 ${\bf CP}^1$ has particular points $0, \infty \in {\bf CP}^1$,
 and let $ 0 \in D(1) \subset S^2= {\bf CP}^1$ be the hemisphere. 
 We choose and fix the following data:
 \vspace{2mm} 
 
 \begin{itemize}
\item a large $l \geq 1$,  

\item a non trivial homotopy class  $\alpha \in  \pi_2(M)$, and

\item  different fixed points
  $p_0, p_{\infty} \in M_0 \subset M $.
  \end{itemize}
  
\vspace{2mm} 
Let $L^2_{l+1}(S^2, M_i)$ be the sets of
$L^2_{l+1}$ maps from $S^2$ to $M_i$.
  Let us introduce the spaces of Sobolev maps:
\begin{align*}
 {\frak B}_i \equiv  {\frak B}_i (\alpha) & = \{  \ u   \in   L^2_{l+1}(S^2, M_i) :  
 [u] = \alpha ,   \\
 & \int_{D(1)} u^*(\omega) = \frac{1}{2}<\omega, \alpha>, \quad
  u(*) = p_* \in M_0  , * \in \{ 0, \infty \}  \ \}.
  \end{align*}

 Let $E(J)_i, \ F_i \mapsto S^2 \times M_i$ be 
 vector bundles whose fibers are respectively:
\begin{align*}
& E(J)_i(z,m) = \{ \phi: T_z S^2 \mapsto T_mM_i : 
   \text{ anti complex linear } \}, \\
& F_i(z,m) = \{ \phi: T_z S^2 \mapsto T_mM_i : \text{ linear } \}.
\end{align*}
Then we have   two  Hilbert bundles over ${\frak B}_i$:
 \begin{align*}
&  {\frak E}_i = L^2_l({\frak B}_i^*(E(J)_i)) = 
\cup_{u \in {\frak B}_i}  \ \{u\} \times L^2_l(u^*(E(J)_i)), \\
&   {\frak F}_i = L^2_l({\frak B}_i^*(F_i)) 
= \cup_{u \in {\frak B}_i}  \ \{u\} \times L^2_l(u^*(F_i)).
\end{align*}
There exist compatible, free and continuous  $S^1$ actions 
on these Hilbert bundles,
which are both induced from the standard action  on ${\C} \subset {\bf CP}^1$.

\begin{remark}
One may regard:
 $$E(J) = \cup_{i \geq 0} \  E(J)_i , \qquad F = \cup_{i \geq 0}  \ F_i$$
  are stratified vector bundles
 over $S^2 \times M $ with   $M=  \cup_{i \geq 0} \ M_i$.
 So 
 their unions
$ {\frak E} = \cup_{i \geq 0} \ {\frak E}_i$ and  $     {\frak F} = \cup_{i \geq 0}  \ {\frak F}_i $
are  stratified   by Hilbert bundles
over $ {\frak B} \equiv  \cup_{i \geq 0} \ {\frak B}_i$,
but both $ {\frak E} $ and  $     {\frak F} $
 are not Hilbert bundles.

Later when we analyze structure of  holomorphic maps into $M$,
we have to take another completions on them.
\end{remark}


 The non linear {\em Cauchy-Riemann operators} and the tangent maps
 are defined respectively  as  sections:
 \begin{align*}
 & \bar{\partial}_{J_i} \in C^{\infty}({\frak E}_i \mapsto {\frak B}_i),  
\quad
  \bar{\partial}_{J_i} (u) = Tu + J_i  \circ Tu \circ \sqrt{-1} , \\
& T \in C^{\infty}({\frak F}_i \mapsto {\frak B}_i),  
\quad  T (u) = Tu, \quad   u \in {\frak B}_i
 \end{align*}
 where $i$ is the complex conjugation on $S^2 = {\bf CP}^1$.

\begin{definition}
$u$ is a {\em holomorphic curve}, if
it  satisfies the equation:
 $$\bar{\partial}_{J_i}(u)=0.$$
\end{definition}
Notice that a holomorphic map $u$ into $M_i$ can be 
regarded as the one into $M_{i+1}$, since 
it also satisfies the equation
$\bar{\partial}_{J_{i+1}}(u)=0$
by compatibility condition.

\vspace{3mm}

Let us   define the moduli space of holomorphic curves by:
$$   {\frak M}(\alpha, M_i,J_i)  = 
    \{ \ u \in C^{\infty}(S^2, M_i) \cap {\frak B}_i(\alpha) 
      : \bar{\partial}_{J_i} (u) =0   \ \}.$$
       
       \begin{remark}
       $(1)$        There is an induced $S^1$ free action on 
        ${\frak M}_i \equiv  {\frak M}(\alpha, M_i,J_i) $,
        if $\alpha \in \pi_2(M)$ is non zero.

$(2)$
$u \in {\frak M}_{i_0}$  implies $u \in {\frak M}_i$
for any $i \geq i_0$, since the embedding
$ {\frak M}_{i_0} \subset  {\frak M}_i$ holds.
\end{remark}

We  say that $J$ is {\em regular}
 at $u \in {\frak M}_{i_0}$,  if the linearizations:
 $$D \bar{\partial}_{J_i}(u): T_u {\frak B}_i \mapsto ({\frak E}_i)_u$$
 are onto 
 for all $i \geq i_0$.
 
 $J$ is regular, if it is regular at any $u \in  {\frak M} \equiv \cup_{i \geq 0} \ {\frak M}_i$.

\vspace{3mm}

The following follows  from the Riemann-Roch and 
 the implicit function theorem:
\begin{proposition}\label{RR}
Let $[(M_i,  \omega_i, J_i)]$ be a regular
almost K\"ahler sequence.
Then the moduli spaces are $S^1$  manifolds with the dimension equality:
$$    
\dim  {\frak M}(\alpha, M_i, J_i)  =  2  <c_1(T^{1,0} M_i) , [u]>   +2 \dim  M_i -1 .$$
 ${\frak M}(\alpha, M_i,J_i )$ is compact, if moreover  $\alpha$ is  minimal.
\end{proposition}
Later on  we omit to denote $\alpha$.

\begin{definition}
The moduli space of holomorphic curves into an  almost K\"ahler sequence is given by:
$$  {\frak M}[(M_i,  \omega_i, J_i)] = \cup_{i \geq 0}  \ 
        {\frak M}(\alpha, M_i,J_i).$$
\end{definition}
 ${\frak M}[(M_i, \omega_i, J_i)]$     is an  $S^1$ 
stratified manifold, if $[(M_i,  \omega_i, J_i)]$ is  a regular
almost K\"ahler sequence.

\begin{example}
Let us consider the standard holomorphic embedding
${\bf CP}^1 \hookrightarrow {\bf CP}^n$
with fixed two points. Modulo $S^1$ action, this is the unique element  
in the moduli space which is regular in the minimal class.
\end{example}

\subsection{Sacks-Uhlenbeck's estimates}
\begin{lemma}\label{hol-est}
Let $[(M_i, \omega_i, J_i)]$ be a minimal almost K\"ahler sequence.
Then  there are constants $c_l \geq 0$ so that
any element $u \in {\frak M} [(M_i, \omega_i, J_i)]$ 
satisfies  the uniform estimates:
$$|\nabla^{l} u|_{C^0(S^2)}   \leq  c_l.$$
\end{lemma}

\begin{proof}
We verify  only the uniform estimate $|\nabla u|_{C^0} \leq c$.
The estimates on the higher devrivatives 
follow from the elliptic regularity.

  There is a biholomorphic isomorphism:
$$\Phi: Z   = {\mathbb R} \times S^1    
 \cong {\bf CP}^1 \backslash  \{ 0, \infty \}, \quad (r,t) 
\to exp(r + 2\pi  \sqrt{-1} t)$$
where we equip with the standard complex structure on $Z$. 
For any holomorphic curve 
$u \in {\frak M}$, let us  regard it as:
 $$u: {\mathbb R} \times S^1 \mapsto M$$ with
$u(- \infty) = p_0$ and $u(\infty) = p_{\infty} \in M_0$.

It follows from  the holomorphic condition
$\frac{\partial }{\partial s} u + J \frac{\partial }{\partial t}u=0$ 
that  the equalities hold:
\begin{align*}
||du||^2 & = \omega(\frac{\partial }{\partial s} u, J \frac{\partial }{\partial s}u)
+ \omega(\frac{\partial }{\partial t} u, J \frac{\partial }{\partial t}u) \\
& =2 \omega(\frac{\partial }{\partial s} u , \frac{\partial }{\partial t}u)
= 2 ||u^*(\omega)||^2.
\end{align*}

\begin{sublemma}[SU]\label{SU}
There  are  constants $C$  and $\epsilon >0$
determined by $[(M_i, \omega_i, J_i)]$  so that 
for any holomorphic disk $u : D^2 \mapsto M = \cup_{ i \geq 0} \  M_i$ and
$E = \int_{D^2} u^*(\omega) \leq \epsilon$, 
the estimate holds:
$$ \varphi (x) \leq C E, \quad \varphi = |du|^2$$
for all $x \in D^2(\frac{1}{2})$.
 \end{sublemma}

{\em Proof of  lemma \ref{hol-est}:}
Let us fix a small positive constant $\delta >0$.
Then for any $u \in {\frak M}([(M_i, \omega_i, J_i)])$, we put
 $s(u) \equiv s_{\infty}(u)  - s_0(u)>0$, where:
\begin{align*}
& s_0(u) =  \sup \{ \ s \in {\mathbb R}:
\ d(\ u(( - \infty , s) \times S^1) ,\  p_0 \ ) \ \leq \ \delta \ \} , \\
& s_{\infty}(u) = \inf  \{ \ s \in {\mathbb R} :\
 d( \ u((s, \infty) \times S^1) , \ p_{\infty} \ ) \ \leq \ \delta \ \}.
 \end{align*}

{\bf Step 1:}
We claim that for $0< \mu \leq    \frac{s(u)}{3}$, 
there is a positive $\epsilon >0$ determined by 
$[(M_i, \omega_i, J_i)]$   and $\mu$ with the estimates:
$$ \int_{ (- \infty, s_0(u)+ \mu] \times S^1} u^*(\omega) ,  \ \ 
 \int_{[s_{\infty}(u)-  \mu, \infty) \times S^1} u^*(\omega) \ \  \geq \epsilon.$$
We  verify the first estimate only. 
The latter follows by  the same argument.

Notice that  the translation on $Z$ is an automorphism
(but it  does not preserve  the required condition 
$\int_{D(1)} u^*(\omega) = \frac{1}{2}<\omega, \alpha>$ on ${\frak B}$).

Let us choose a translation $T$ on $Z$ 
so that $s_0(u \circ T) = 0$ holds. Notice $s(u \circ T) =s(u) \geq   3 \mu$.
Then one may assume $s_0(u)=0$, since  the equality: 
$$ \int_{ (- \infty, s_0(u \circ T)+ \mu] \times S^1} (u \circ T)^*(\omega) =
 \int_{ (- \infty, s_0(u)+ \mu] \times S^1} u^*(\omega)$$
holds.
Let $D^2(b) \subset S^2$ be the disk with the radius $b >0$.
Then we choose   $  a >0$ as:
 $$ (- \infty, s_0(u)+ \mu] \times S^1 =  D^2(1+a) 
\backslash 0  \subset S^2$$ where we identify 
 $ (- \infty, s_0(u)] \times S^1 =  D^2(1) \backslash \{0\}$. 
 We put  $D=D^2(1)$ and $D' = \ D^2(1+a)$.
 
   Let us put $B_{\delta}(0) \equiv \{ m \in M: d(p_0, m) < \delta) \} \subset M$
   as $\delta$ neighborhood of $p_0$.
   Then 
   $u(s,t) \in  \partial   B_{\delta}(0) $ and so 
$d(u(s,t), u(-\infty)) = \delta $ holds
    at $s= s_0(u)$ and some $t \in S^1$
   with $(s,t) \in \partial D$.

Suppose $\int_{D'} u^*(\omega) < \epsilon$ 
could hold for 
sufficiently  small $\epsilon  =\epsilon(\mu)>0$.
 Then by sublemma \ref{SU}, 
  the uniform estimate on the derivative: 
 $$ |du|   \leq C(\mu) \sqrt{\epsilon} $$ should hold at any 
 point of $D$. 
 This is a contradiction if $\epsilon >0$ is too small, since  $u(s,t) \in \partial B_{\delta}(0)$ 
 and $d(p_0, u(s,t)) = \delta$ as above.
 
This verifies  the claim.

\vspace{3mm}

{\bf Step 2:}
Let us  proceed  by contradiction argument.
So suppose contrary. Then there are families
$\{ u_i\}_i \subset {\frak M}[(M_i, \omega_i, J_i)]$
and $\{p_i\}_i \subset S^2$ with 
$ |\nabla u_i|(p_i) \to \infty$.
As \cite{hofer and viterbo}  page $611$, one may assume the two properties:
\begin{align*}
&  |\nabla u_i|(x) \ \leq \ 2 |\nabla u_i|(p_i), \\
& \lim_{i \to \infty} \ \epsilon_i |\nabla u_i|(p_i) = \infty
\end{align*}
for all $x$ with  $d(x, p_i) \leq \epsilon_i $, where 
$\epsilon_i \to 0$  as $i \to \infty$.

Let $D_i = D_i(p_i)$ be small  balls with the center $p_i$.
 We rescale these  balls  as
 $B_i = |\nabla u_i|(p_i) D_i(p_i)$
 by multiplying  $|\nabla u_i|(p_i)$ on their radii,
 where one regards $B_i \subset {\C}$. 
 By conformal invariance, one gets a family of 
holomorphic maps $v_i : B_i \mapsto M = \cup_{i \geq 0} \ M_i$.
This family satisfies uniform bounds:
 $$  |d v_i | (p_i) =1,     \quad 
|d v_i|(x)    \leq 2 \quad \text{ for }      |x-p_i| \leq   \epsilon_i |\nabla u_i|(p_i)   .$$
In particular by choosing  small $1 >> a, \epsilon' >0$, 
 $$| |dv_i|^2(x)  - |dv_i|^2(p_i) | \leq  \epsilon'$$
hold for all $x \in D_a(p_i)$ by  elliptic regularity,
where $D_a(p_i) \subset B_i$ is $a$ ball with the center $p_i$,
and $a$ is independent of $i$.
 This gives   the lower bounds
$|dv_i|(x) \geq \sqrt{1 -  \epsilon'}$.
So  the uniform estimates hold from below:
$$\int_{D_{b_i}(p_i)} v_i^*(\omega) \geq \int_{D_a(p_i)} v_i^*(\omega) \geq C >0$$ 
for all $b_i \geq a$ with $D_{b_i}(p_i) \subset B_i$.

{\bf Step 3:}
On the other hand
 apriori bounds 
$\int_{B_i} v_i^*(\omega) \leq m$ hold from above
where  $m$ is the minimal invariant.
We claim that there is 
 some family 
$ R_i  \leq  \epsilon_i   |\nabla u_i|(p_i)$ with  $R_i \to \infty$
 such that the length $\delta_i$ of
$x_i \equiv v_i(R_i \exp(2 \pi i t)): S^1  \mapsto M$
must decay   $\delta_i \to 0$.

In fact since $[a,b_i] \times S^1 \subset B_i$ hold with $b_i \to \infty$,
  there are some $R_i$ so that
the decay:
$$\int_{[R_i-1, R_i+1] \times S^1} v_i^*(\omega) \to 0$$
must hold. Then  the decay
$\sup_{x \in R_i \times S^1} |dv_i|(x) \to 0$ holds
by sublemma \ref{SU}, which verifies the claim.

{\bf Step 4:}
Thus  there is a family of small disks $\{ d_i \}_i \subset M$
which span $x_i$, and $\int_{d_i} \omega \to 0$.
Let $B_i' \subset B_i$ be $R_i$ balls with the center $p_i$,
whose  boundaries  are $x_i$.
Let us put two `almost'  holomorphic
spheres:
$$  u_i' =  
\begin{cases}
 & u_i
  \text{ on } S^2  \backslash B_i'  \\
  &   d_i
    \end{cases} ,
    \quad \quad 
  v_i' = B_i' \cup d_i.$$
By the condition, these must satisfy:
\begin{align*}
& <\omega, u_i'> + <\omega, v_i'> \to m >0, \\
& \lim_i <\omega, u_i'> , \quad 
 \lim_i <\omega, v_i'> \ \  \geq \ 0
\end{align*}
By minimality, one of $<\omega, u_i'>$ or $<\omega, v_i'>$ 
must be zero for all large $i$.
By step $2$ and  $3$,
  $<\omega, v_i'>$ must be positive and equal to $m$. So
  $<\omega, u_i'> =0$ must hold.

First of all, suppose  there is a uniform lower bound
$s(u_i) \geq 3\mu_0 >0$.
There are three cases;

(1) an infinite subset of $\{ p_i\}_i$ is contained in $(-\infty , s_0(u_i)] \times S^1$
or

(2)  is contained in  $[s_{\infty}(u_i), \infty) \times S^1$ or

(3) in $[s_0(u_i), s_{\infty}(u_i)] \times S^1$.
\vspace{3mm} \\
Suppose the case (1). Then by step $1$,
there is a positive $\epsilon >0$
with
$ \int_{[s_{\infty}(u_i) - \mu_0, \infty) \times S^1} (u_i')^*(\omega) \geq \epsilon$.
This implies the asymptotic bounds:
$$\lim_i <\omega, u_i'> \  \geq \  \epsilon$$
which give a contradiction.
The other cases can be considered similarly.

\vspace{3mm} 

{\bf Step 5:}
Let us verify that 
  $s(u_i) \to 0$ cannot happen.
  This will complete the proof of the lemma. 
Suppose contrary.
Let us take $p =0, q=\frac{1}{2} \in S^1$.
Then since $u_i(o \times s_*(u_i) ) \in B_{\delta}(*)$,
 $* = 0, \infty$ and $o = p,q$,
and since
 $d(B_{\delta}(0), B_{\delta}(\infty))  >0$ is positive,
 there are families
 $\{t_i \}$ and $\{r_i \}$, $t_i, r_i  \in [s_0(u_i), s_{\infty}(u_i)] $,
such that $|\nabla u_i|(p \times t_i ), |\nabla u_i|(q \times r_i) \to \infty$.
On the other hand one has a lower bound
 $d(p \times t_i, q \times r_i) \geq \frac{1}{2}$
in ${\mathbb R} \times S^1$.
By the same arguments as step $2,3,4$, 
one obtains two non trivial almost holomorphic spheres,
which also cannot happen by  minimality of the homotopy class. 

\end{proof}

\subsection{Hilbert completion of function spaces}\label{Hilbert-compl}
Let us  introduce the basic function spaces 
on the infinite dimensional analysis.

Recall the Sobolev $l+1$ space $ {\frak B}_i(\alpha)$
of maps from sphere to $M_i$ in \ref{f.d.pre}.
Let   us take an element:
$$u \in {\frak B}_i(\alpha) \ \subset \ 
 {\frak B}(\alpha)  \ \equiv \ \cup_{i \geq 1} \  {\frak B}_i(\alpha)$$
 and let $U(u) \subset  {\frak B}(\alpha)$
 be  a small  neighborhood of $u$ in the set of $L^2_{l+1}$
maps from $S^2$ to $M$.
Below we will describe its completion to a Hilbert manifold $\hat{U}(u)$.

Let us check the Sobolev embedding of maps into Hilbert space.
\begin{lemma}
There is a constant $c_l$ with the uniform estimate:
$$||u||_{C^{l-1}(S^2) } \leq c_l ||u||_{L^2_{l+1}(S^2)}.$$
\end{lemma}

\begin{proof}
By uniformity of complete local charts, 
it is enough to verify 
 the uniform estimate:
$$||u||_{C^{l-1}_c}  \leq c_l ||u| |_{L^2_{l+1}}$$
for $u \in C_c( D^2 ; H ) $ with the open unit disc $D^2 \subset {\mathbb R}^2$.

The Sobolev estimate
$||\tilde{u}||_{C^{l-1}_c(D^2)} \leq c_l ||\tilde{u}||_{L^2_{l+1}(D^2)}$ holds 
for scalar functions $\tilde{u} \in C_c(D^2)$.
Let $H $ be the closure of ${\mathbb R}^{\infty}$ with the standard norm,
and express $u =(\tilde{u}_0, \tilde{u}_1, \dots) \in C_c^{l-1}(D^2; H)$.
Then we have the estimate:
\begin{align*}
\sum_{k=0}^{l-1} \
|\nabla^k u|^2 (m) & = \sum_{k=0}^{l-1} \
 \sum_{j \geq 0} \ |\nabla^k \tilde{u}_j|^2(m) \\
& \leq c_l  \sum_{j \geq 0} \ || \tilde{u}_j||^2_{L^2_{l+1}(D^2)} = c_l ||u||^2_{L^2_{l+1}(D^2)}
\end{align*}
for any $m \in D^2$. By taking supremum of the values in the left hand side,
we obtain the desired estimate.
\end{proof}

\begin{remark}
These Hilbert Sobolev spaces
 admit the
free and continuous $S^1$ actions. In precise there is a constant $C>0$
with the inequalities:
$$ C^{-1} ||u|| \ \leq \ \sup_{t \in S^1} ||tu|| \ \leq \ C ||u||$$ for all elements $u$ 
in such  a space.
\end{remark}

\vspace{3mm}

Let us precisely describe how to equip with the Sobolev norm
 on ${\frak B}(\alpha)$.
Let  $ \varphi(p) :D (\epsilon) \equiv \cup_i \ D^{2i}(\epsilon)
     \hookrightarrow  M = \cup_{i \geq 0} \  M_i$ be a
  complete almost K\"ahler chart at $p$.
Sometimes we will identify $D(\epsilon)$ with  $D(p)$ where: 
$$D(\epsilon) 
\ \subset \ {\mathbb R}^{\infty} \ \subset  \ H , \quad
D(p) \equiv   \varphi(p)( D(\epsilon) )  \ \subset \ M.$$
$D(\epsilon) $ admits  the induced metric from $ \varphi(p)$,
which is 
uniformly equivalent to the standard one on $H$.

Let us fix  the following data;
\vspace{2mm} 

$(1)$  finite set
of points $s_0, \dots, s_k \in S^2$, 
\vspace{2mm}

$(2)$ an open cover $U_0, \dots, U_k$ with
$s_i \in U_i \subset S^2$, and 
\vspace{2mm}

$(3)$ a partition of unity $f_0, \dots, f_k$
 over $S^2$.
 \vspace{3mm}

  For  each $u \in {\frak B}(\alpha)$, 
one can choose  large $k$ so that 
each  image  $u(U_j)$ is contained in
a complete almost K\"ahler chart at $\varphi(p_j)$ with $p_j = u(s_j)$.
Then one can express its restriction as:
$$u|U_j : (U_j,s_j) \mapsto (D(p_j), p_j).$$
Identifying $D(p_j)$ with  $D(\epsilon)$ as above, one may regard these maps as:
$$u|U_j: (U_j,s_j) \mapsto (D(\epsilon), 0) \subset (H,0).$$
Let  $u \in U(u) \subset {\frak B}(\alpha)$
be a small  open subset.  
Then locally any element $v \in U(u) $
can be expressed as $v|U_j : U_j \mapsto {\mathbb R}^{\infty} \subset H$.
We introduce the Sobolev norm on $ U(u)$ by:
$$||v||^2_{L^2_{l+1} } = \sum_{0 \leq j \leq k} \sum_{0 \leq a \leq l+1}   \
        \int_{U_j} |\nabla^a ( f_j v)|^2 (m) dm.$$
By completion, one obtains  the Hilbert manifold
$ \hat{U}(u) $
 which contains   a neighborhood  of $u  \in {\frak B}(\alpha)$.
 Then we put the Hilbert manifold:
 $$\hat{\frak B}(\alpha)  = \cup_{ u \in {\frak B}(\alpha)}  \ \hat{U}(u).$$
 
Notice that if $u$ is holomorphic, then $k$ above can 
be chosen uniformly by lemma \ref{hol-est}.

In a similar way, let us  introduce
the Hilbert norm on the set of sections of $u^*(E(J))$ as
follows;
let us take any $\varphi \in \Gamma(u^*(E(J)))$. 
Then one can express
 the restriction  as:
$$\varphi|U_j:  TU_j \mapsto  TD (\epsilon)
= D (\epsilon) \times {\mathbb R}^{\infty}$$
which is anti linear with respect to $(i, J_{u( m)})$ at $(m,u(m))$.
Notice that $\varphi$ is of the form
$\varphi(m) = (u(m), \varphi_m)$ with $\varphi_m  \in E(J)(m,u(m))$.
Then  define:
$$||\varphi||^2_{L^2_l} \equiv  \sum_{0 \leq j \leq k}  \sum_{ 0 \leq a \leq l} \
\int_{U_j}  |\nabla^a(f_j \varphi|U_j)|^2 (m) dm.$$

Let $W(u, \varphi) \subset {\frak E}$ be an open neighborhood of $(u, \varphi)$.
By taking completion with respect to the above norm, 
one obtains the Hilbert bundle:
$$ \hat{W}(u, \varphi) \to \hat{U}(u).$$
Then we put the total Hilbert bundle:
$$\hat{\frak E} (\alpha) \equiv  \cup_{(u,\varphi) \in {\frak E}} \  \hat{W}(u, \varphi)
 \ \to \ \hat{\frak B}(\alpha).$$
Similarly we obtain the Hilbert bundle:
$$\hat{\frak F} (\alpha) 
 \ \to \ \hat{\frak B}(\alpha).$$

\vspace{3mm}

Now the Cauchy-Riemann operator extends to a smooth section as:
$$\bar{\partial}_J  : \hat{\frak B}(\alpha)  \to  \hat{\frak E} (\alpha)$$
which restricts as:
$$\bar{\partial}_J  : \hat{U}(u)  \to   \hat{W}(u, \partial_J(u)) .$$

Let us denote: 
$$ \hat{\frak M} [(M_i,  \omega_i, J_i)] = 
\cup_{u \in {\frak M} [(M_i, \omega_i, J_i)] }
   \ \  \{ \ v \in  \hat{U}(u)
      : \bar{\partial}_{J} (v) =0  \ \}.$$
 This space is apriori bigger than the moduli space
${\frak M}[(M_i,\omega_i,J_i)]$.
Nonetheless later we verify  their coincidence each other 
under some conditions.

\subsubsection{Some functional analysis}
Let $H$ be a Hilbert space and $L \subset H$ be a closed linear subspace.

\begin{lemma}\label{func-analy}
Let $F: H \to H$ be a bounded operator with closed range,
whose kernel consists of finite dimensional subspace.
Then $F(L) \subset H$ is also closed.
In particular  $F(L)$ is closed if $F$ is injective.
\end{lemma}

\begin{proof}
If kernel $F=0$, then $F: H \cong F(H)$ gives an isomorphism.
In particular $F(L)$ is closed.

Suppose ker$(F) = K \subset H$ is of finite dimension.
Then $F$ induces an isomorphism
$F: H/ K \cong F(H)$, where we equip with the metric 
on $H/K$ by use of orthogonal decomposition $H = K^{\perp} \oplus K$.
Then it is enough to see that the image of the projection 
$pr(L) \subset H/K$ is still closed.

One may assume that $L \cap K= 0$ by replacing $L$ by $(L \cap K)^{\perp}$  in $L$,
when it  has positive dimension.

Suppose a sequence $\{ \bar{v}_i\}_i \subset pr(L)$ converge to some element
$\bar{v} \in H/K$. 
By the assumption, their representatives $v_i \in L$ of $\bar{v}_i$ are unique.
Let us represent  $v_i  = v_i^1+v_i^2 \in L$ with 
 respect to the decomposition $H=K^{\perp}\oplus K$.

We claim that $||v_i||$ are uniformly bounded. 
Suppose contrary and assume $||v_i|| \to \infty$.
Then by normalizing as $w_i = ||v_i||^{-1}v_i = w_i^1+w_i^2$,
both convergence $||w_i^1|| \to 0$ and $||w_i^2|| \to 1$ should hold.
Since $K$ is finite dimensional and $L$ is closed,
a subsequence $w_i$ converges to some element $w \in L \cap K$ with $||w||=1$.
This contradicts to  our assumption, which verifies the claim.

Now since $\{v_i^2\}_i \subset K$ is a bounded sequence,
a subsequence converges to some element $v^2 \in K$.
Since $v_i^1$ converges to $v$, it follows from these 
that a subsequence of  $\{v_i\}_i$ converges to $v+ v^2 \in L$.
This implies $\bar{v} \in pr(L)$.
\end{proof}

\begin{remark}
The assumption of finite dimensionality is necessary.
Let $H$  be a separable  infinite dimensional Hilbert space,
and choose an orthonormal basis $\{v_i\}_i$.
Let $0< a_i  \to 0$
be a decreasing family of numbers.

Let us consider a surjective bounded map:
$$F = \text{ id } \oplus 0 : H \oplus H \to H$$
and a closed subspace $L $
spanned by the basis:
$$L = \text{ span } \{ w_i = (a_i v_i, v_i): i =0,1,2, \dots\} \subset   H \oplus H.$$
We claim that the image of the restriction $F|L$ is not closed.
Suppose  contrary. Then since $F|L$ is injective, the restriction must be
an isomorphism by the open mapping theorem.
So there must exist some $C>0$ with the uniform estimates:
$$|a_iv_i| = |F(w_i)| \geq C|w_i| = C \sqrt{a_i^2+1}. $$
But the left hand side converge to $0$,
which cannot happen.
\end{remark}

The following abstract property is a key to our Fredholm theory 
we develop later:
\begin{corollary}
Suppose the above situation, and choose another Hilbert space $W$.
Then the image of the Hilbert space tensor product $L \otimes W$
of  the induced operator $F \otimes 1: H \otimes W \to H \otimes W$,
 still has closed range.

In particular if $F$ is an isomorphism, then  $F \otimes 1$ is also the same.
\end{corollary}

\begin{proof}
 Let us put $E= L \cap \ker  (F)$, and decompose 
 $L \cong L' \oplus E$. Then  $F(L)=F(L')$ holds.
 Since the restriction $F|L'$ is injective, 
 it gives the isomorphism onto $F(L)$
 by the open mapping theorem.
 
Since  the restriction $F \otimes 1|L' \otimes W$ gives the isomorphism 
onto $F \otimes 1(L' \otimes W)= F\otimes 1(L \otimes W)$,
the conclusion follows.
\end{proof}

\subsection{Geometric conditions}\label{geo-cond}
We study  functional analytic properties of the Cauchy-Riemann operators 
over almost K\"ahler  sequences which satisfy 
the  geometric conditions we have introduced in \ref{sym}.

Our aim in \ref{geo-cond} is to verify the following:

 \begin{theorem}\label{main-thm}
 Let $[(M_i, \omega_i, J_i)]$ be a symmetric  K\"ahler sequence.
 
$ (1)$
Suppose it is  regular, and $\dim \ \cup_{i \geq 0}    \ker  D_u \bar{\partial}_i =N$ is finite.
 Then it is in fact strongly regular of index $N$.

 In particular  
  $ {\frak M}[(M_i, \omega_i , J_i)] $ is a regular $N$ dimensional  manifold.
  
 $ (2)$ If   moreover  it is  isotropic
 and  is regular with respect to a minimal class  $\alpha \in \pi_2(M)$, then 
   the equality holds:
   $$ {\frak M}[(M_i, \omega_i , J_i)]=  {\frak M}( M_0,\omega_0, J_0) .$$
   Moreover it is  compact.
     \end{theorem}
  
 \begin{proof}
 $(1)$  follows from combination of lemma \ref{ker-cond} with proposition \ref{s-reg.cond}.

$(2)$ is verified in lemma \ref{moduli-struc}.
\end{proof}

\subsubsection{Strong regularity over symmetric K\"ahler  sequences}
Let $[(M_i, \omega_i, J_i)]$ be a symmetric
  almost K\"ahler sequence, and 
  choose its symmetric  data  $\{ (P_i, \pi_k)\}_{i,k}$ with respect to $(M_k, M_{k+1})$
 in \ref{sym}.

For any $u \in {\frak B}_k \subset {\frak B}_{k+1} \subset {\frak B}$,
let $\hat{U}(u) \subset {\frak B}$ be as in \ref{Hilbert-compl}.
Let us put
  $U(u)_{l}  = {\frak B}_{l} \cap \hat{U}(u)$.
  
  There is the  extended projection:
  $$\bar{\pi}_k : \hat{U}(u) \to U(u)_k$$
  with $\bar{\pi}_k|U(u)_k =$ id, given by the composition
  $v \to \pi_k \circ v$.
  Then
 the isomorphism:
 $$T_u U(u)_{k+1} \cong  T_u {\frak B}_k \oplus V(u)_k$$
holds, where:
$$V(u)_k= \ker  (\bar{\pi}_k)_*  \cap T_uU(u)_{k+1} .$$

\begin{lemma}
The complete isomorphism holds:
$$T_u \hat{U}(u) \ \cong  \ T_u {\frak B}_k \ \oplus  \ V(u)_k \otimes H$$
where $H$ is a  separable Hilbert space.
\end{lemma}

\begin{proof}
This follows from  the symmetric property and lemma \ref{exp}.
 \end{proof}

The Cauchy-Riemann operator $\bar{\partial}_J$ and  the tangent map $T$
 give smooth sections respectively:
$$  \bar{\partial}_J : \hat{U}(u) \mapsto \hat{\frak E}|\hat{U}(u), \quad
  T:  \hat{U}(u) \mapsto  \hat{\frak F}|\hat{U}(u).$$

\begin{definition} 
Let  $[(M_i,  \omega_i, J_i)]$ be
 a regular  almost K\"ahler sequence.
It is strongly regular,  if the differential:
$$D\bar{\partial}_u: 
T_u \hat{U}(u) \mapsto T_u  \hat{ \frak E} $$ is transverse to the $0$
section $\hat{U}(u) \subset \hat{ \frak E}$
 for  any $ u \in {\frak M}[(M_i, \omega_i, J_i)]$.
\end{definition}

\begin{lemma}\label{closed-range}  
Let $[(M_i, \omega_i, J_i)]$ be a symmetric 
   K\"ahler sequence. 
   
   Then
  $D\bar{\partial}_J  : T_u \hat{U}(u) \mapsto  T_u\hat{\frak E}$ 
 has  closed range. 
\end{lemma}

\begin{proof}

{\bf Step 1:}
 Let us take
$u \in   {\frak M}(M_k, \omega_k, J_k)$, and 
 $P_i:(M, M_k) \cong (M,M_k)$ be the 
symmetry data
 for all $i \geq  k+1$.
 
 Recall the notations in \ref{f.d.pre}, and 
consider the bundle over $S^2 \times M_k$:
 $$F^{\perp}_{k}(m,z)= \{ \phi: T_z S^2 \to \text{ Ker }( \pi_k )_*
 \cap T_m M_{k+1}: \text{ linear } \} \subset F_{k+1}(z,m).$$
  Then we put the Hilbert sub bundle over ${\frak B}_k$:
 $${\frak F}^{\perp}_{k}= L^2_l ({\frak B}_k^*(F^{\perp}_{k}))= \cup_{u \in {\frak B}_k} \{u\} \times L^2_l(u^*(F^{\perp}_{k})) \subset {\frak F}_{k+1}| {\frak B}_k.$$
 
 There is a   bundle decomposition:     
   $${\frak F}_{k+1} |{\frak B}_k \cong {\frak F}_k \oplus {\frak F}_{k}^{\perp}$$
  over ${\frak B}_k$ given by:
  $$\bar{\phi} \to ( \ (\pi_k)_*(\bar{\phi}), \ \bar{\phi} -  (\pi_k)_*(\bar{\phi}) \ ).$$
 It follows from  symmetric property that  
 the bundle decomposes:
 $$\hat{\frak F}|{\frak B}_k \ \cong \ {\frak F}_k \ \oplus \ {\frak F}_{k}^{\perp} \otimes H$$ 
   as lemma \ref{decomp}.
 Let:
   \begin{align*}
   DT_{k+1}  =DT_k \oplus DT_{k}^{\perp} 
    : T_u {\frak B}_{k+1} & = T_u {\frak B}_k \oplus  V(u)_k \\
&  \to T_u {\frak F}_{k+1} = T_u {\frak F}_k \oplus  
 T_u {\frak F}_{k}^{\perp}
 \end{align*}
  be the tangent map.
Then the  total tangent map is described as:
 \begin{align*}
 DT & = DT_k  \oplus (DT_{k}^{\perp} \otimes \text{ id }) \\
&  : T_u \hat{U}(u) 
   \ \cong  \ T_u {\frak B}_k  \ \oplus \ V(u)_k \otimes H \\
& \  \  \to \   T_u \hat{\frak F} \ \cong \ T_u {\frak F}_k \  
\oplus  \ T_u {\frak F}_{k}^{\perp}  \otimes H.
\end{align*}

{\bf Step 2:} 
Let us consider $J$.
Notice that it always preserves the restriction
as $V(u)_k \to V(u)_k$, since it commutes with $(\pi_k)_*$. 
However 
it is not always the case that its restriction gives a self map as
$T_u {\frak F}_k^{\perp} \to  T_u {\frak F}_k^{\perp}$.
On the other hand  it happens if  it is integrable.
In fact  the formula:
$$ D\bar{\partial}_J(v) = (DT + J \circ DT \circ \sqrt{-1})(v) +  N(v) $$
holds, where $N$ involves $\nabla J$
and  $ N \equiv 0$ when $J$ is integrable.
So
$   D\bar{\partial}_J   = DT + J \circ DT \circ i  $ holds  if it  is K\"ahler.

In particular if we decompose these function spaces
by use of the holomorphic local charts as in step $1$,
then $D\bar{\partial}_J $ can be also expressed 
as a form $K_1 \oplus (K_2 \otimes \text{ id })$.

Now consider the composition with the projection:
 $$D\bar{\partial}_l
 : T_u {\frak B}_l \to T_u {\frak F}_l \to  ({\frak F}_l)_u.$$ 
 This is Fredholm for any $l \geq k$,
 which follows from the well known analysis of holomorphic curves 
 into finite dimensional symplectic manifolds (see [HV]).
In particular  the map
 $D\bar{\partial}_l
 : T_u {\frak B}_l \to T_u {\frak F}_l$
  has closed range with finite dimensional kernel.

 Since 
 $V(u)_k \subset T_u{\frak B}_{k+1}$ is a closed linear subspace, 
it follows from   lemma \ref{func-analy} that 
 $K_2 \otimes \text{ id }$ has closed range.
So the direct sum
 $K_1 \oplus (K_2 \otimes \text{ id })$ also has closed range.
\end{proof}

 \subsubsection{Index computations}
Let $[(M_i, \omega_i, J_i)]$ be an almost K\"ahler sequence,
and $\pi_k: U_{\epsilon}(M_k) \to M_k$ 
be the  holomorphic projection
with $\pi_k|M_k=$ id from 
 a  small neigborhood  in $ M= \cup_{i \geq 0} \  M_i$ for each $k \geq 0$.
For  $u \in {\frak B}_k$, 
let $\bar{\pi}_j: \hat{U}(u) \mapsto {\frak B}_j$ be 
the  induced   projections for all $j \geq k$.

Let us compare two operators: 
$$\bar{\partial}_J : \hat{U}(u) \mapsto \hat{\frak E}|\hat{U}(u), \qquad
\bar{\partial}_i  :      {\frak B}_i   \mapsto {\frak E}_i .$$ 

\begin{lemma}\label{ker-cond} 
Let
 $[(M_i, \omega_i, J_i)]$ be an almost K\"ahler sequence.

If    $\cup_{i \geq 0}   \ker D_u\bar{\partial}_i$ is of finite dimension,
   then the equality holds:
$$   \ker D_u\bar{\partial}_J    \ =  \ \cup_{i \geq 0} \    \ker  D_u\bar{\partial}_i.$$

In particular the left hand side is of finite dimension.
\end{lemma}

\begin{proof}
The condition implies    
$\cup_{i \geq 0}  \ker D_u\bar{\partial}_i =  \ker D_u\bar{\partial}_{i_0}$ 
for some $i_0$.

Suppose contrary and assume  $\ker D_u\bar{\partial}_J    
\ne \cup_{i \geq 0}   \ker D_u\bar{\partial}_i$.
Let $u_t \subset \hat{U}(u)$ be a smooth curve with $u_0 =u$ and 
$u_t'|_{t=0} \equiv v \in  \ker  D_u\bar{\partial}_J  $
but $v \not\in  \cup_{i\geq 0}   \ker D_u\bar{\partial}_i$.
  It follows from the equality:
 $$D_u\bar{\partial}_j(\bar{\pi}_j(v)) = \bar{\pi}_j(D_u\bar{\partial}(v))  =0$$ 
  for  all $j \geq k$
that  $(\pi_j)_*(v) $ lies in $\ker D_u \bar{\partial}_j$. 
Hence it  must be contained in  $\ker D_u\bar{\partial}_{i_0}$.
 Since $j$ is arbitrary, this implies $v \in  \ker D_u\bar{\partial}_{i_0}$. 
 This is a contradiction.
 \end{proof}
 \vspace{3mm}

\begin{proposition}\label{s-reg.cond}
Let $[(M_i, \omega_i, J_i)]$ be
 a symmetric  K\"ahler sequence.
  Let us  choose any $u \in {\frak M}(M_k, \omega_k, J_k)$.
  
If the uniform bound  $\dim \coker  D_u\bar{\partial}_i \leq M$ 
holds  for any $i\geq k$,
then  $\dim \coker  D_u \bar{\partial}_J \leq M$ also holds.

In particular if it is regular, then it is in fact strongly regular.
\end{proposition}

\begin{proof}
$D\bar{\partial}_J$ has closed range by  lemma \ref{closed-range}.
Suppose  $\dim \coker D_u\bar{\partial}_J \geq M+1$ could hold, 
and take orthonormal  elements  $u_1, \dots, u_{M+1}$
in $\coker D_u\bar{\partial}_J$.
There is a  large $l >>k$ so that  
$u_i^l = \pi_l \circ u_i$ is  defined for any $1 \leq i \leq M+1$.

For  small $\epsilon >0$, let us choose sufficiently large  $l$ so that  the estimates below
hold, where $B \subset \text{ im}D_u \bar{\partial}_l \ \subset  ({\frak E}_l)_u$ 
is the unit ball:
$$||u_i^l||^2 \geq 1 - \epsilon, \quad| < u_i^l, u_j^l>| \leq \epsilon,
\quad  |< B , u_i^l>| \leq \epsilon.$$
There are numbers $a_1, \dots, a_{M+1} \in {\mathbb R}$ with 
  $\sum_{i=1}^{M+1} |a_i|^2 =1$
such that  $v \equiv \sum_{i=1}^{M+1} a_i u_i^l $ lies in im $ D_u\bar{\partial}_l$, 
since  $\dim \coker D_u\bar{\partial}_l \leq M$ holds.
Let us pick up $i$ with $|a_i| = \sup_{1 \leq j \leq M+1} |a_j| \geq \frac{1}{\sqrt{M+1}}$.
Then one should have the estimates:
$$\epsilon \geq |<v,   u_i^l>| \geq |a_i| (1- \epsilon) -  \epsilon \sum_{i \ne j} |a_j|
\geq |a_i|(1- \epsilon) - \sqrt{M} \epsilon.$$
Since $\epsilon$ can be  arbitrarily small, this  is a contradiction.
\end{proof}

\begin{example}
  ${\bf CP}^{\infty}$  is strongly regular of index $1$ by proposition \ref{RR}.
  \end{example}

So for a regular and symmetric  K\"ahler sequence, 
the moduli space of  holomorphic curves is strongly regular
 with the expected index.

 Notice that the strong regularity condition is   stable under small perturbations,
while just regularity is not the case in general.
      With  the above analysis, we would like to propose the following:

\begin{conjecture}
 Let $[(M_i, \omega_i, J_i)]$ be 
a symmetric K\"ahler sequence.

$(1) $
 One can perturbe the complex structure (to be almost K\"ahler)
 so that the result  could become strongly regular.
 
 $(2)$  
  index $ D\bar{\partial}_J = M$ holds  when
 index $D\bar{\partial}_i =M $  for all large  $i$  and 
  $\cup_{i \geq 0}  \ker D\bar{\partial}_i$  is of finite dimension.
\end{conjecture}

\subsubsection{Compactness of moduli spaces}

   \begin{lemma}\label{moduli-struc}
 Let  $[(M_i,  \omega_i, J_i)]$ be an isotropic   K\"ahler sequence,
 which is regular with respect to a minimal class  $\alpha \in \pi_2(M)$.
Suppose moreover
$\cup_{i \geq 0}    \ker D_u\bar{\partial}_{l_0} < \infty$
is of finite dimension.

 Then  the equality holds:
 $$ {\frak M} [(M_i,  \omega_i, J_i)]    = {\frak M}( M_0,  \omega_0, J_0) .$$
 Moreover they are  compact.
 \end{lemma}
 

  \begin{proof}  
 Let us choose an element $[u] \in  {\frak M} [(M_i,  \omega_i, J_i)] $.

   {\bf Step 1:}
Let us verify that  there is some $l_0$ so that
    the connected component ${\frak M}(u)$
  containing $u$ has all their images  in $M_{l_0}$.
  
  It follows from 
  proposition \ref{s-reg.cond} 
  that the moduli space is strongly regular.
  By lemma \ref{ker-cond}, 
  there is some $l_0$ so that  the equality
  $   \ker D_u\bar{\partial}_J    = \cup_{i \geq 0}    \ker D_u\bar{\partial}_{l_0}$
  holds.
  
  Suppose there is some $u' \in {\frak M}(u)$ whose image is not contained in $M_{l_0}$.
  Then take a smooth path $u_t$ between $u$ and $u'$ in ${\frak M}(u)$ for $t \in [0,1]$.
  There should exist $ t_0 \in [0,1]$ such that the image of $u_{t_0}$ lies in $M_{l_0}$,
  but it is not the case for any  $u_t$ with $t \in (t_0, t_0+\epsilon]$,
  where $\epsilon >0$ is a positive number. 
  Choose some $t \in (t_0, t_0+\epsilon]$ and $l_0' \geq l_0$ such that 
  the image of $u_t$ lies in $M_{l_0'}$.
  
  Now $[u] \in {\frak M}( M_{l_0},  \omega_{l_0}, J_{l_0})$ and 
  $[u_t] \in {\frak M}( M_{l_0'},  \omega_{l_0'}, J_{l_0'})$,
  which are both regular  manifolds of the same dimension $N$.
  At $[u]$, the tangent spaces coincide:
  $$T_{[u]}  {\frak M}( M_{l_0},  \omega_{l_0}, J_{l_0})
  = T_{[u]}  {\frak M}( M_{l_0'},  \omega_{l_0'}, J_{l_0'})$$
  with the inclusion:
  $${\frak M}( M_{l_0},  \omega_{l_0}, J_{l_0})
\subset {\frak M}( M_{l_0'},  \omega_{l_0'}, J_{l_0'}).$$
So the local charts at $[u]$ for these moduli spaces must coincide.
On the other hand by lemma \ref{hol-est},
 one may assume that 
  $[u_t]$ lies in the local chart of ${\frak M}( M_{l_0'},  \omega_{l_0'}, J_{l_0'})$ at $[u]$.
  This cannot happen.

  {\bf Step 2:}
Suppose  there could exist  some   $k$   such that  all elements in 
  ${\frak M}(u)$  have  their images  in $M_{k+1}$
   but some element does not have its image   in $M_k$.
  
  Let   $P_i^t: (M, W_i,M_k) \cong (M, M_{k+1},M_k)$ be the isotropies for $i \geq k+1$,
  where $M=\cup_{i \geq 0}  M_i$. There is some $u' \in {\frak M}(u)$ so that
   the image  $P^0_i(u') =u'$ is contained in $M_{k+1}$, 
  while $(P^1_i)^{-1}(u')$ are not the case for all  $i \geq k+2$.
  This implies that the image of  ${\frak M}(u)$ cannot be contained in $M_{k+1}$,
  since $(P^1_i)^{-1}(u')$  also consists of an element   in ${\frak M}(u)$.
  This  contradicts to the assumption.
  So ${\frak M}(u)$ must be contained in $M_k$.
  
  Next let us replace the pair $(k,k+1)$ by $(k-1, k)$.
  Then  the image of  ${\frak M}(u)$ is contained in $M_{k-1}$
  by the same argument.
    
  Let us continue this process. Then finally we find that the image of
  ${\frak M}(u)$ must be contained in $M_0$.
  \end{proof}

It would be interesting to study more general case 
when  the equality:
 $$ {\frak M} [(M_i,  \omega_i, J_i)]    = {\frak M}( M_k,  \omega_k, J_k) $$
 holds for some $k$, if 
 $ {\frak M} [(M_i,  \omega_i, J_i)]   $ is a smooth manifold of finite dimension.

\section{Hamiltonian dynamics}

\subsection{Bounded Hamiltonians}
Let  $[(M_i, \omega_i, J_i)]$ be an almost K\"ahler sequence.
 A bounded function $f:M  \to {\mathbb R}$   is called a  {\em bounded Hamiltonian}, if it is 
of completely bounded geometry
on  $M = \cup_{i \geq 0} M_i$ (see \ref{var-norm}).

Let 
$\{ \ \varphi(p): D(\epsilon)\hookrightarrow \cup_{i \geq 0}  M_i \ \}_{p\in M}$
be  a uniformly bounded covering
by $\epsilon$ complete almost K\"ahler charts.
By pulling back the bounded Hamiltonian
as $\varphi(p)^*(f):  D(\epsilon) \to {\mathbb R}$,
let us regard the restriction of the differential $df$ as  a one form on $D(\epsilon)$.
The Hamiltonian vector field $X_f$ on $D(\epsilon)$ is defined as the unique
vector field  which obeys the equality:
$$-df(Y) = \omega(X_f, Y)$$
for any  vector field $Y$ of completely bounded geometry.

Below we  introduce three classes of   bounded Hamiltonian functions on 
$M= \cup_{i \geq 0} \ M_i$.

\vspace{2mm}

\begin{definition}\label{ham.fcn-def}
$(1)$  A bounded Hamiltonian function
$f: M \to [0, \infty)$ is {\em pre-admissible}, if there are open  neighborhoods 
 $N_0, N_{\infty } \subset M$ of $p_0, p_{\infty} \in M_0$ respectively
so that 
$f|N_0 \equiv 0$ and $f | N_{\infty} \equiv \sup f$
hold.

\vspace{2mm}

$(2)$
 $f$ is {\em proper}, if for any $i$, there is 
$j$ so that: 
$$(df)_m \in T^*_m M_j$$ holds 
for any $m \in N(M_i) \subset M$, where $N(M_i)$ is 
a  neighborhood of $M_i$ in $M$.

Moreover  the following holds:
$$ \lim_{k \to \infty} ||f - f_k||_{C^{\alpha}(M_k) }=0$$
for any $\alpha \geq 0$,
where $f_k \equiv f|M_k$ are the restrictions.
Notice that:
 $$ \lim_{k \to \infty} F|M_k  = F_k$$ also holds in $C^{\alpha-1}$, 
where $F_k: M_k \cong M_k$ are the Hamiltonian diffeomorphisms
with respect to the restrictions $f_k$.

\vspace{2mm}

$(3)$
 $f$ is  $l_0$-{\em connected}, if
there is an open and connected subset  $p_0 \in U \subset M$ which contains
$\supp d f \subset U $ such that its closure $\bar{U} \subset M$ is $l_0$-connected:
$$\pi_l(\bar{U} \cap M_k) =0 \quad (l \leq l_0)$$
for all $k =0,1, \dots$,
where $\supp d f $ is the closure of $\{ m \in M : df(m) \ne 0\}$.

\end{definition}

\vspace{2mm}

\begin{lemma} Let $[(M_i, \omega_i, J_i)]$ be an almost K\"ahler sequence, and 
take a proper bounded Hamiltonian 
$f: M= \cup_{i \geq 0}  M_i \to {\mathbb R}$.

 Then there is the parametrized diffeomorphisms as its integral:
 $$F_t : M \cong M$$
 which preserve the symplectic form.
 \end{lemma}

 Notice that $M$ is not complete, since it is  countable union 
 of finite dimensional manifolds and  modeled on ${\mathbb R}^{\infty}$.
 So  without properness,
 $f$ cannot induce maps on $M$ to itself in general.

 \begin{proof}
 Let  $f: M \to [0, \infty)$ be a  pre-admissible
  bounded Hamiltonian.   
   Take any $m \in M_i$, and 
 consider the gradient vector field $X= \text{grad }f$ on $N(M_i) \subset M$.
  Then the restriction $X|_{M_i} $ takes the values in $TM_j$.
  In fact any vector $W \in (T_m M_j)^{\perp} \subset T_m M$
 is orthogonal to $X$, since the equality holds:
 $$<X, W> = df(W).$$
 So the restriction of the Hamiltonian vector field $X_f = - J \circ \text{grad }f$ 
 on $M_i$ also takes 
 the values in $TM_j$.
 It follows from uniqueness of the integral
 that range of the Hamiltonian diffeomorphisms $F_t(m)$ land on $M_j$ 
 for any $m \in M_i$ and all small $0 \leq t \leq t_0$.
 The conclusion follows immediately.
 \end{proof}

 \vspace{2mm}

We call $F=F_1$  as the {\em Hamiltonian diffeomorphism}.
Notice that  $F_t(U) \subset U$ hold
if supp $d f \subset U \subset M$.

\begin{example}
Let $D \subset {\mathbb R}^{\infty}$ 
and $D^k  \subset {\mathbb R}^k$ 
be the unit disks with $D= \cup_{k \geq 0}D^k$, and consider the bounded Hamiltonian
$f:D \to {\mathbb R}$ given by:
$$f(x_0,x_1, \dots) = \sum_{k=0}^{\infty}  \ x_kx_{k+1}.$$
$f$ is  not  proper.
 
Let $\rho: {\mathbb R} \to {\mathbb R}$ be  a smooth function
such that $\rho(a) \equiv 0$ for $|a| \leq \epsilon $ and $\rho(a) =a$
 for $|a| \geq 2 \epsilon$ for some $\epsilon >0$.
Then $h$ below is proper:
$$h(x_0,x_1, \dots) = \sum_{k=0}^{\infty}  \ x_k \rho_k(x_{k+1}).$$
\end{example}

The topological condition $(3)$ above is used when we introduce the invariant of the action functional
in \ref{action-func}.

\subsection{Symplectic action functional}\label{action-func}
Let $M$ be a finite dimensional manifold, and 
take a connected open subset $p_0 \in U \subset M$ with a fixed point.

A {\em cone} of $(p_0,U)$ is a smooth map
$u:(-\infty, 0] \times S^1 \to M$
so that $u(-\infty, S^1) = p_0$ and  $u(0, S^1) \subset U$ hold.

\begin{lemma}
Suppose  
  $U$ is $2$-connected.
A cone of  $(p_0, U)$ canonically defines an element in $\pi_2(M)$.
\end{lemma}

\begin{proof}
Let $\Omega (M)$ be the based loop space with the base $p_0$.
$u(s, \quad): S^1 \to M$ gives a path in the free loop space, and 
by use of the path $u(\quad, 0)$, it is  lifted to a path 
$\tilde{u}: (- \infty, 0] \to \Omega(M)$.

Let us put $u(0,0)=m \in U$ and choose another path
$\gamma $ from $m$ to  $p_0$  in $U$.
The loop $u(0, \quad): S^1 \to U$ can be contracted since $U$ is simply connected.
So one can obtain an extension $\tilde{u}$  so that 
$\tilde{u}: (- \infty, 1] \to \Omega(M)$ is the family of loops along $\Psi \equiv u( \quad , 0) \cup \gamma$,
and   $\tilde{u}(1)$ is the union $\Psi \cup - \Psi$ which is the trivial element.
So  $\tilde{u}(- \infty) = \tilde{u} (1)\in \Omega(M)$ holds,
and $\tilde{u}$ gives an element $\alpha \in \pi_1(\Omega(M))= \pi_2(M)$.

$\alpha$ is independent of choice of contraction of $u(0, \quad)$ since $\pi_2(U)=0$ holds.
It is also independent of choice of $\gamma$ since $\pi_1(U)=1$.
\end{proof}

\vspace{3mm}

Let us call the triplet
$(\alpha, p_0,U)$
as the {\em action functional data},
where $\alpha \in \pi_2(M)$ is the element defined as above from a cone of $(p_0, U)$.

\vspace{3mm}

Let $(M, \omega)$ be a  finite dimensional symplectic manifold, and choose
a pre-admissible Hamiltonian
$f: M \to [0, \infty)$ which is $2$-connected  over
$p_0 \in  U \supset \supp df$ with $f(p_0)=0$.
Let $F_t : (M, U)  \cong (M, U)$ be the Hamiltonian path with $F=F_1$,
and consider a fixed point $m \in U$ by $F$.
$F_t$ fixes  a neighborhood of $p_0$ by pre-admissibility.

 Since $m$ is fixed by $F$,
  there exists 
 a non trivial periodic orbit $u(0, \quad) : S^1 \to U$
 of the Hamiltonian vector field $X_f$ with $u(0,0)=m$.
Then  take 
a cone:
$$u:(-\infty, 0] \times S^1 \to M$$
with $u(-\infty, S^1) = p_0$ and with  the induced element $\alpha$.
The triplet $(\alpha, p_0,U)$ gives an  action functional data.

\subsubsection{Small perturbations}\label{small-perturb}
Before going further, let us consider small perturbations.
For $\epsilon >0$, let $\mathcal{F}_{\epsilon}(f) $ 
be the set of smooth functions on $M$ such that the estimates:
$$||f-g||_{C^1} < \epsilon$$
hold for  $g \in \mathcal{F}_{\epsilon}(f) $.

\begin{lemma}\label{perturb} 
Let $U$ be as above for $f$ which is $2$-connected.
For $\delta >0$, there is $\epsilon >0$ such that 
any periodic orbit $l$ with:
 $$ \text{ length } l \geq \delta$$
of Hamiltonian diffeomorphisms $G_t: M \cong M$ 
is contained in $U$ 
for any $g \in \mathcal{F}_{\epsilon}(f) $.
\end{lemma}

\begin{proof}
Notice that there is $\mu >0$ such that 
the upper bound $\mu \geq $ length $l$ holds, since
it satisfies the equation
$\frac{dl}{dt} = X_G(l(t))$ and $C^0$ norm of the vector field 
is uniformly bounded for any $g \in \mathcal{F}_{\epsilon}(f) $. 

If $\epsilon >0$ is sufficiently small, then $l$ must be contained in $U$,
since the $C^0$ norm of $X_G$ on $U^c$  is smaller than $\epsilon$.
\end{proof}

\subsubsection{Action functional}\label{action-functional}
Let $u$ be a cone as above.
Consider the curve $\gamma : [0,1] \to M$  by $u(\quad , 0)$, 
and  the loop $l$ given by two arcs $ - \gamma$ and $F \circ \gamma$.
Let us consider   the disc $\Delta $ with $\partial \Delta = l \cup - u(0, \quad)$ by:
$$\Delta : (- \infty, 0] \times [0,1] \to M, \quad
\Delta(s,t) = F_t(u(s,0))$$
Then with respect to the data $(\alpha, p_0, U)$, we define the action functional  by:
$$\delta(F;m) = \int_{u} \omega +  \int_{\Delta} \omega.$$
Notice that the homological boundary of $u +\Delta$ is $l$, and so 
the equality:
$$\delta(F;m) = \int_l \mu$$
holds by Stokes theorem, where $d \mu =\omega$ on $U$.

\vspace{3mm}

The proofs of the following two lemmas extend the arguments in [P].

\begin{lemma}
If $U$ is $2$-connected, then
the action functional is independent of choice of  cones
and  Hamiltonian paths
with respect to the data $(\alpha, p_0,U)$.
\end{lemma}
So it makes sense to denote $\delta(F;m)$ when $U$ is $2$-connected,
 once the data $(\alpha, p_0, U)$ is given.

\begin{proof}
{\bf Step 1:}
Let us verify that it is  independent of choice of  cones.
Let us choose another cone $u'$ such that $u'(0, \quad)=u(0, \quad)$
 is the periodic orbit 
with respect to $X_f$. The loop consisted by $u(\quad ,0)$ with $u'(\quad ,0)$
can be spanned by a disk $D$ in $M$, 
since 
 we have fixed the element $\alpha \in \pi_2(M)$. Moreover 
$u'$ can be deformed to $u$ rel  $(-\infty , S^1) \cup (0, S^1)$.
Then we consider the sphere:
$$S = u \cup  \Delta \cup -u' \cup - \Delta' \cup D \cup -F(D).$$
 $[S] =0 \in \pi_2(M)$ vanishes since
both $\Delta'$ and $u'$ can be deformed 
to $\Delta$ and $u$ respectively, passing through $D$ and $F(D)$.

Then  the equalities:
$$\int_{u} \omega +  \int_{\Delta} \omega
- ( \int_{u'} \omega +  \int_{\Delta'} 
\omega) =- \int_D \omega +  \int_{F(D)} \omega=0$$
hold where the last one  holds since $F$ preserves the symplectic form.

So $\delta(F;m)$ is independent of choice of  cones
with fixed  data $(\alpha , p_0, U)$.

{\bf Step 2:} 
Let us verify that it is
 independent of choice of  Hamiltonian paths.
Let us take another Hamiltonian path $F_t'$ with $F_1' =F$,
and consider the loop $F'_t(m)$.
There is a map $v: [0,1] \times S^1 \to U$ with $v(0, t) =F_t(m)$ and $v(1,t)=F'_t(m)$,
since $U$ is simply connected.
Then let us choose another 
cone $u'$ which is obtained by concatenation of $u$ with $v$.
Clearly the corresponding element is the same as $u$ in $\pi_2(M)$.
We may assume the equality  $u=u'$ on $(-\infty, -\epsilon] \times S^1$,
and $u,u' : [- \epsilon ,0] \times S^1 \to U$
for some $\epsilon >0$.
Let us put: 
$$\Sigma = u \cup  \Delta, \ \Sigma' = u' \cup  \Delta' \  \subset  \ M.$$
$\Sigma \cup - \Sigma'$ consistes of a sphere $S$.
Let us put:
$$K = \{ (s,t) \in (- \infty, 0] \times [0,1]: \Delta(s,t) \ne \Delta'(s,t), \}.$$
Then $ \Delta (K) \cup \Delta'(K)$
 is contained in $U$,
since the support of the differential of the Hamiltonian functions 
are contained in $U$.

We use a general fact that for an $n$ dimensional  CW pair $(A,B)$ with $n-1$-connected $Y$,
any continuous map $h: B \to Y$ can be extended over $A$.
Let us 
 apply the above fact to the pair
 $A= K \times [0,1]$ with $B= \partial K  \times [0,1] \cup K \times \{0,1\}$
 with $n=3$ and:
  $$h \equiv \Delta'|K \times \{0\} \cup \partial K  \times [0,1], \ \Delta| K \times \{1\} \to U.$$
  This gives a continuous deformation from $S$ to $\Sigma \cup - \Sigma$
  which is contractible.
  So we obtain the equality:
  $$\int_{\Sigma} \omega = \int_{\Sigma'} \omega.$$
 \end{proof}

\vspace{3mm}

Let $f: M \to [0, \infty) $ be a pre-admissible Hamiltonian function with a periodic orbit 
 $\{ F_t(m)\}_{t \in [0,1]}$  in $U$. Then 
we  define the symplectic action:
$${\frak A}(f, m) = \int_u \omega - \int_0^1 f(F_t(m))dt .$$

Let us fix the data $(n \alpha, p_0, U)$ for $F^n$ and $n =1,2, \dots$

\begin{lemma}\label{linearity}

$(1)$
 $\delta(F; m) = {\frak A} (f,m) = \int_{u} \omega - \int_0^1 f(F_t(m))dt$.

$(2)$
 $\delta(F^n ; m) = n \delta(F; m)$ for all $n=1,2, \dots$

$(3)$
 Suppose $U$ is $2$-connected.
If  $\delta(F;m) \ne 0$, then
$||dF^n||_{C^0(U)}$ must grow at least linearly  with respect to  $n$.
In particular 
 the cyclic group generated by $F$ is infinite.
\end{lemma}

\begin{proof}
See proposition $2.4.A$ in [P] for $(1)$.
Notice that in our case $p_0$ is fixed under Hamiltonian deformations 
by $F_t$.

Let us verify $(2)$. 
It is enough to see $n=2$. Let us put:
$$\Delta^2 : (- \infty, 0] \times [0,2] \to M, \quad
\Delta^2(s,t) = F_t(u(s,0))$$
and  $u^2 \equiv u(s,2t)$. 
$\Delta^2$ is just $\Delta \cup F(\Delta)$ where
$\Delta : (- \infty, 0] \times [0,1] \to M$ with 
$\Delta(s,t) = F_t(u(s,0))$.
In particular:
$$u^2 \cup  \Delta^2 = (u \cup \Delta) \cup (u \cup  F(\Delta))$$
Since $F$ preserves the symplectic form, the equalities hold:
$$\int_{u \cup  F(\Delta)} \omega =
\int_{F(u \cup \Delta)} \omega = \int_{u \cup  \Delta} \omega.$$
So we obtain the equalities:
$$\delta(F^2, m) = \int_{u^2 \cup  \Delta^2} \omega = 2 \int_{u \cup  \Delta} \omega = 2 \delta(F,m).$$

For  $(3)$, $H^2(U : {\mathbb R})=0$ by Hurewicz isomorphism theorem.
Let $\mu_0$ be a primitive one form of $\omega$ with $d \mu_0 = \omega$ on $U$, 
and put
$C_0 = ||\mu_0||_{C^0(U)} < \infty$.

Let $l_n$ be  the loops in $M$ consisted by 
$- \gamma$ and $F^n(\gamma)$ with $\gamma =u(\quad, 0)$.
Notice the estimates:
$$\text{ length } l_n \leq l(\gamma)(1 + ||dF^n||_{C^0 (U)})$$
where $l(\gamma)$ is the length of $\gamma$.
Then we have the inequalities:
$$n | \delta(F; m)| = | \int_{ l_n}  \ \mu_0|
\leq C_0 \text{ length } l_n \leq 
C_0 l(\gamma)(1 + ||dF^n||_{C^0(U)}).$$
\end{proof}

\subsubsection{Perturbed case}
Recall \ref{small-perturb} and take a  bounded Hamiltonian 
$g$ which is sufficiently close to $f$ in $C^{\infty}$.
Let $G: M \cong M$ be the Hamiltonian diffeomorphism
with respect to $g$.
Let $l$ be the periodic orbit as in lemma \ref{perturb}.
Take a cone $u$ with $p_0$.

There is a slightly different situation to \ref{action-functional},
since $p_0$ is not fixed by $G$ and the point moves a little.
We define the same action functional 
$\delta(G,m)$.

\begin{corollary}
The lower bound:
$$n \delta(G,m) \leq C_0 l'(\gamma)(1 + ||dG^n||_{C^0 (U)})$$
holds for each $n \geq 1$.
\end{corollary}

\begin{proof}
Notice that the homological boundary in this case is
$l'$ which is given by union of three arcs,
$- \gamma$, $G_t(p_0)$ and $F \circ \gamma$.

The same argument tells us the equality:
$$\delta(G^n, m) = n \delta(G,m).$$

Let $l_n'$ be the loop given by the union of $-\gamma$ with $G^n (\gamma)$.
Then we have  the estimates:
$$\text{ length } l'_n \leq l'(\gamma)(1 + ||dG^n||_{C^0 (U)}).$$

Combining these two inequalities, we obtain the lower bound:
$$n \delta(G,m) \leq C_0 l'(\gamma)(1 + ||dG^n||_{C^0 (U)}).$$
\end{proof}

\begin{remark}\label{act-func.esrtimate}
Let $\varphi: M \to \R$ be a cut off function with 
$\varphi|N_{\frac{\epsilon}{2}}(p_0 \cup p_{\infty}) \equiv 0$
and $\varphi|N_{\epsilon}(p_0 \cup p_{\infty} )^c \equiv 1$. Then 
$g'  \equiv \varphi g $ 
is also a small perturbation of $f$,
and is pre-admissible.
Moreover the estimate holds:
$$n \delta(G',m) \leq \frac{n}{2} \delta(G,m).$$

Suppose $l$ is a periodic orbit of $G_t$.
It follows from the proof of lemma \ref{perturb} that  the intersection
is empty:
$$l \cap N_{\epsilon}(p_0) = \phi.$$
So $i$ is also a periodic orbit of $G'$.
\end{remark}

\subsection{Proof of theorem \ref{unif.cyc}}
Let  $[(M_i, \omega_i, J_i)]$ be an
 almost K\"ahler sequence, and
 $f:M  \to [0, \infty)$  be a non constant proper bounded Hamiltonian
 with $F : M \cong M$. Suppose $f$ is
 pre-admissible with fixed points 
 $p_0,p_{\infty} \in M_0\subset M = \cup_{i \geq 0} \  M_i$.

Let $f_i = f|M_i$ be the restrictions and $F_i : M_i \cong M_i$
be the associated Hamiltonian diffeomorphisms.
Our key lemma is given by the following:

\begin{lemma}\label{key lemma}
Assume the conditions in theorem \ref{unif.cyc}.

Then there are $\epsilon >0$,  uniformly bounded positive numbers:
 $$C \geq \lambda_i \geq c >0$$ 
 from both sides, and a family of  fixed points $m_i \in U$ with respect to $\tilde{F}_i$
such that the uniform bounds hold from below for all sufficiently large $i >0$:
$$|\delta(\tilde{F}_i ;  m_i)| > \epsilon$$
where $\tilde{F}_i$ are the Hamiltonian diffeomorphisms with respect to $\lambda_i f_i$.
\end{lemma}
The proof of the lemma occupies section \ref{cobordism}.
Before going into the proof, let us verify theorem  \ref{unif.cyc} assuming the key lemma 
\ref{key lemma}.

\vspace{3mm} 

\begin{proof}\label{proof of unif.cyc}
Let us verify  theorem \ref{unif.cyc}.
Notice that  $F_i$ is
 different from $F|M_i$  in general, since it is
 Hamiltonian with respect to $f|M_i$, while
  the image $F(M_i)$ does not necessarily coincide with $M_i$.

{\bf Step 1:}
Firstly we 
 verify that there is some $\epsilon ' >0$ 
so that:
  $$ ||d(F^m)||_{C^0(M_i)} \geq  m \epsilon'$$ hold
for any $m $ and $i \geq i_0(m)$.
This implies that $F$ is infinitely cyclic.

Let $\lambda_i$ be in lemma \ref{key lemma} and 
$\lambda = \liminf_i   \lambda_i  \geq c >0$ be the positive number,
where one  takes a subsequence if necessary.
Let $F_t$ and $\tilde{F}_t$ be the Hamiltonian diffeomorphisms
with respect to $f$ and $\lambda f$ respectively.
Notice the relation 
$\tilde{F}_t = F_{\lambda t}$.

By properness, 
 $\{ \tilde{F}_i \}_i$ converges to $\tilde{F}$ in $C^1$, 
 where $\tilde{F}_i$ are in lemma \ref{key lemma}.

{\bf Step 2:}
There are fixed points 
$m_i \in M_i$ with respect 
to $\tilde{F}_i :M_i \cong M_i$ so that 
the uniform estimates:
$$|\delta(\tilde{F}_i ; m_i)| > \epsilon$$
hold for all large $i >0$ by lemma \ref{key lemma}.

Let $\gamma_i$ be  paths between $p_0$ and $m_i$ in $U$,
and consider the  loops $l_n^i$ consisted by $\gamma_i$ with $(\tilde{F}_i)^n(\gamma_i)$ as in lemma \ref{linearity}.
Since  diameter of $M$ is  finite,
one may assume that the lengths of $\gamma_i$ are uniformly bounded from above. 

Then as in the proof of lemma \ref{linearity}$(3)$, we have the uniform estimates:
$$n \epsilon \leq C_0  (1+ ||d(\tilde{F}_i)^n||_{C^0(M_i)}).$$

{\bf Step 3:}
 For any $l_0$:
 $$\lim_{i \to \infty} ||F^{l_0}|M_i- (F_i)^{l_0}||_{C^1(M_i)} =0$$
 hold  since  $f$ is proper.
 
 Suppose $\lambda_i =1$.
Then there are  $p_i \in M_i$ so that 
the estimates:
 $$n \epsilon \leq C_0 |d(F_i)^n|(p_i)$$
  hold for all large $n \geq 0$.
Choose $n_0$ so that 
$||d(F_i)^{n_0l_0}||_{C^0(M_i)} \geq 2$ 
hold for all large $i$.

Take  large $i >>1$ so that the estimates
$||F^{n_0l_0}|M_i - (F_i)^{n_0l_0}||_{C^1(M_i)} <1$ hold.

Then   
we get a contradiction, if $F^{l_0}$ could be the identity.

{\bf Step 4:}
Next let us verify the general case. Let us consider
 $ \tilde{F}_i$  in step $2$.
 Then 
  the estimates:
$$n \epsilon \leq C_0 |d(\tilde{F}_i)^n|(p_i) \qquad (*)$$
  hold 
  for all large $n \geq 0$ as in step $3$.
 
 There are families of numbers 
 $0 \leq \mu_n < 1 $ and $a_n \in {\mathbb N}$
 such that  $n \lambda_i = a_n +\mu_n$.
   Since  the equalities
  $(\tilde{F}_i)^n \equiv \tilde{F}_i \circ \dots \circ \tilde{F}_i
  =(\tilde{F}_i)_{n} = (F_i)_{\lambda_i n}$
  hold, we obtain the estimates  for all $n \geq 0$:
  \begin{align*}
 a_n \epsilon   & \leq   \lambda_i n \epsilon \leq C_0 \lambda_i |d(\tilde{F}_i)^n|(p_i) \\
& =
  C_0 \lambda_i |d(F_i)_{\lambda_i n}|(p_i) 
  =  C_0 \lambda_i |d(F_i)_{\mu_n} \circ d(F_i)^{a_n}|(p_i) 
  \leq C_0' \lambda |d(F_i)^{a_n}|(p_i) .
  \end{align*}
So we obtain the estimates
$a_n \epsilon \leq C_0'' |d(F_i)^{a_n}|(p_i)$ for some $C_0''$.
Since $a_n \to \infty$ as $n \to \infty$, 
we can repeat the same argument as step $3$ above.

{\bf Step 5:}
Choose a small $\mu >0$ and take any $g \in \mathcal{F}_{\mu}(f)$.
Notice that $g$ may not be pre-admissible, but is $\mu$-close in norm.
So by remark \ref{act-func.esrtimate}, a similar estimate as $(*)$ in step $4$
also holds for the Hamiltonian diffeomorphisms $G_i$ of $g_i$,
by replacing $n$ by $2n$ on the left hand side, if necessary.
The rest argument is the same, and obtained infinite cyclicity of $G$.
\end{proof}

\vspace{3mm} 

\begin{remark}\label{low.bound}
In fact   the uniform bound:
$$\lambda_i \ \leq \ <u, \omega> (\sup f - \inf f)^{-1}$$ holds
by  lemma \ref{top.bound}.
\end{remark}

\section{Periodic orbit and cobordism}\label{cobordism}
Let  $[(M_i, \omega_i, J_i)]$  be an almost K\"ahler sequence, 
and put $M = \cup_{i \geq 0} \ M_i$.
 Let us  fix the following data;
 (1) a large $l \geq 1$,  
(2) a non trivial homotopy class  $\alpha \in  \pi_2(M)$ and
(3)  two different  points
  $p_0, p_{\infty} \in M_0 \subset M \equiv  \cup_{i \geq 0} \  M_i$.

  Let us  take
a pre-admissible Hamiltonian
$f: M \to [0, \infty)$ 
with $f|N_0 \equiv 0$ and $ f|N_{\infty} \equiv \sup f$. 
There is $\delta >0$ so that the open neighborhoods 
$N_0, N_{\infty} \subset M$ both contain $\delta$ balls with the centers
$p_0,p_{\infty}$ respectively.
We denote its restrictions by 
$f_i: M_i \to {\mathbb R}$.

 ${\bf CP}^1$ has particular points $0, \infty \in {\bf CP}^1$,
 and let $ 0 \in D(1) \subset S^2= {\bf CP}^1$ be the hemisphere.

\subsection{Finite dimensional case}\label{f.d.case}
Let us recall the construction of the functional
 ${\frak F}_i : {\frak B}_i \to {\frak E}_i$ 
in [HV].
Let us put:
$$W_i = (S^2 \backslash \{0, \infty\} ) \times M_i  \ \cup \
 \{0\} \times  \ (N_0 \cap M_i) 
\ \cup \ \{\infty\} \times \ (N_{\infty} \cap M_i).$$
Notice that $(z,u(z)) \in W_i$ for $u \in {\frak B}_i$.
Let us define the complex anti-linear map:
$$\hat{f}_i(z,m) : T_zS^2 \to T_m M_i \quad (z,m) \in W_i$$
by (1) $\hat{f}_i(0, \quad) =0$, $\hat{f}_i(\infty, \quad) =0$, and
(2) $ \hat{f}_i(z,m)(z) =  \frac{1}{2\pi}\nabla f_i(m)$, where we use the 
identity chart over ${\mathbb C} \subset S^2$.
$\hat{f}_i$ are uniquely determined by anti-linearity.
 Then we define the functional:
 $${\frak F}_i(u) (z) = \hat{f}_i(z,u(z)).$$

Let us consider the cobordism:
$${\frak C}_i = \{  \ (\lambda, u)  \in [0, \infty) \times {\frak B}_i : 
\bar{\partial}_i(u) + \lambda {\frak F}_i(u) =0   \ \}$$
which contains
 ${\frak M}(M_i, \omega_i, J_i)$ by embedding $u \to (0, u)$.
Let us put:
\begin{align*}
& {\frak C}_i (\lambda)= \{  \ u \in {\frak B}_i :  (\lambda, u) \in {\frak C}_i  \ \}, \\
& {\frak C}_i(0 \leq \lambda <\epsilon) = 
\cup_{0 \leq \lambda <\epsilon} \ {\frak C}_i(\lambda).
\end{align*}

There is a bi-holomorphic isomorphism:
$$\Phi: Z   = {\mathbb R} \times S^1    
 \cong {\bf CP}^1 \backslash  \{ 0, \infty \}, \quad (r,t) 
\to exp(r + 2\pi  it)$$
where we equip $Z$ with  the standard complex structure.
Then any $u \in {\frak B}$ can be regarded as a map:
 $$u: {\mathbb R} \times S^1 \mapsto M$$ with
$u(- \infty) = p_0$ and $u(\infty) = p_{\infty} \in M_0$.

Let us put:
\begin{align*}
 & s_0(u) =  \sup \{  \ s \in {\mathbb R}:
\ u(( - \infty , s) \times S^1)  \subset N_0 \  \} , \\
& s_{\infty}(u) = \inf  \{ \ s \in {\mathbb R} :\
 u((s, \infty) \times S^1) \subset N_{\infty} \  \}
 \end{align*}
with  $s(u) \equiv s_{\infty}(u)  - s_0(u)>0$.
Then it admits a continuous $S^1$ action which is induced  from $S^1$ coordinate action
on $\R \times S^1$.

\begin{lemma}\label{top.bound}
Suppose $ {\frak M}(M_i, \omega_i, J_i)$ is  compact, regular and $S^1$
 freely cobordant to non zero
with respect to a minimal class $\alpha$.
Then the followings hold:

$(1)$ ${\frak C}_i$ is  non compact, and
${\frak C}_i(\lambda)$ is empty for: 
$$\  \sup f - \inf f  \  >   \  \lambda^{-1}<\omega_i, \alpha>  .$$

$(2) $
Non trivial periodic solutions $x_i$ exist with respect to $\lambda_i f_i$ 
for some $\lambda_i > 0$, where they are obtained as:
$$x_i \ = \  \lim_{l \to \infty} \ u_i^l(s^l_i, \quad) $$
 for  some divergent sequences
 $u_i^l \in {\frak C}_i(\lambda_i^l)$ with respect to $l$. Here 
 $\lambda_i = \lim_l \lambda_i^l$
 and 
 $\lim_l s(u_i^l) = \infty$ with  $s_0(u_i^l) \leq s_i^l \leq s_{\infty} (u_i^l)$.
\end{lemma}
These are verified in proposition $2.6$, $2.7$  and page $618$ in [HV].

\subsection{Sacks-Uhlenbeck type estimates}\label{SU-type}
In \ref{SU-type}, we fix a minimal element.

We verify the following:
\begin{lemma}\label{su-type lem}
There is positive $\epsilon >0$ independent of $i$
so that for any $u \in {\frak C}_i$ 
 the uniform estimates hold:
 $$\int_{(- \infty, s_0(u)]} \omega ,
  \quad \int_{[s_{\infty}(u), \infty)} \omega \ > \ \epsilon.$$
\end{lemma}
The proof of lemma \ref{su-type lem} uses the following:

\begin{proposition}\label{C1-bound}
There is a constant $C_0\geq 0$ so that 
for any $ u \in {\frak C}_i$, 
 the uniform bound:
$$||du||_{C^0({\mathbb R} \times S^1)} \leq C_0 < \infty$$
holds independently of $i$.
\end{proposition}

Before going into the proof, let us finish the proof 
of lemma \ref{key lemma} assuming lemma \ref{su-type lem}
and proposition \ref{C1-bound}.

\vspace{3mm}

{\bf Proof of lemma  \ref{key lemma}:}
Let us  verify uniform positivity
$\delta(\tilde{F}_i; m_i) \geq  \epsilon$, where $(\tilde{F}_i)_t$ correspond to $\lambda_i f_i$
and $m_i =x_i(0)$ in lemma \ref{top.bound}.

For $ u \in {\frak C}_i(\lambda)$, let us regard it as
$u : {\mathbb R} \times S^1 \to M_i$, and put:
$$a(s) = \int_{(- \infty, s] \times S^1} u^*(\omega) - \int_0^1 \lambda f_i(u(s,t))dt.$$
Then:
 $$\frac{da}{ds} = \int_{S^1}
|J(u) \frac{\partial}{\partial t}u +   \lambda \nabla f_i|^2 dt$$ holds, and so $a(s)$ 
is monotone increasing.

Let $(\tilde{F}_i)_t: M_i \cong M_i$ be the Hamiltonian diffeomorphisms with respect to $\lambda_i f_i$, 
and  consider $u_i^l$ and $x_i$ in lemma \ref{top.bound}, where
 $\{u_i^l(s_i^l, \quad)\}_l$ converge to the periodic orbit  $x_i = (\tilde{F}_i)_t(m_i)$
 in $M_i$.

Because the equalities:
$$a(s_0(u_i^l)) = \int_{(- \infty, s_0(u_i^l)]} \omega$$
hold, it follows from monotonicity and lemma \ref{su-type lem}
 that uniform positivity holds:
$$\int_{(- \infty, s_i^l] \times S^1} (u_i^l)^*(\omega) - \int_0^1 \lambda_i^l f_i(u_i^l(s_i^l,t))dt
 \  \geq  \  \epsilon.$$

It follows from lemma  \ref{su-type lem} and 
 proposition \ref{C1-bound} that for a  large $l_0$,
one may modify:
 $$u_i^{l_0} | (- \infty, s_i^l] \times S^1$$
on small neighborhoods of $\{s_i^l \} \times S^1$ 
so that  they consist of the cones:
$$u_i^{l_0}(s_i^{l_0}, t) =x_i(t)=(\tilde{F}_i)_t(m_i)$$ 
 in \ref{action-func}
 with the uniform bounds:
$${\frak A}(\lambda_i f_i, m_i)
= \int_{(- \infty, s_i^{l_0}] \times S^1} (u_i^{l_0})^*(\omega) - \int_0^1 \lambda_i f_i((\tilde{F}_i)_t(m_i))dt
 \ \geq \  \frac{ \epsilon}{2}.$$
Since 
$\delta(\tilde{F}_i;  m_i) ={\frak A}(\lambda_i f_i, m_i) $
hold by lemma \ref{linearity}$(1)$, this verifies uniform positivity.

Next let us verify uniform bounds  $C \geq \lambda_i \geq c >0$ from both sides
in lemma \ref{key lemma}. 
Uniform bound from above follows by lemma  \ref{top.bound}.

Let us verify uniform lower bound.
Suppose contrary and choose a degenerating sequence
$\lambda_i \to 0$, by taking a subsequence if necessarily.
Notice  $a(\infty) = <\omega, u> + \lambda( \sup f - \inf f)$ 
which is uniformly bounded. Moreover $a(s)$ is monotone increasing.
So 
 there is  a family $u_i \in {\frak C}_i(\lambda_i')$
with $s(u_i) \to \infty$, and $x_i = u_i(s_i, \quad) $ satisfy: 
$$ \int_{S^1} |\frac{d x_i}{d t}|^2 dt \to 0$$
for some $s_0(u_i) \leq s_i \leq s_{\infty}(u_i)$ as $i \to \infty$.
In particular diameters of $x_i$ go to zero. 

So
  one can cut $u_i$ along $x_i$, put  small discs on 
   the 
boundary circles, and  obtain two spheres
$v_i^1$ and $v_i^2$ with 
$p_0 \in v_i^1$ and $p_{\infty} \in v_i^2$.

Because $a(s)$ above is monotone increasing
and $\lambda_i \to 0$, it follows from  lemma  \ref{su-type lem} that
both must satisfy uniform positivity:
$$<v_i^1 , \omega> , \  <v_i^2, \omega > \ \  \geq \ \  \frac{\epsilon}{2} .$$
 On the other hand the convergence:
$$\lim_{i \to \infty} <v_i^1 , \omega>  +  <v_i^2, \omega > = <u_i , \omega> $$
hold, and the right hand side is minimal.  This is a contradiction, and
 we are done.
 This completes the proof of lemma  \ref{key lemma}.

\begin{remark}\label{unif.const}
The proofs below verifies that 
both constants $C_0$ in proposition \ref{C1-bound}
 and $\epsilon$ in lemma \ref{su-type lem}
depend only on $||f||_{C^{l+1}(M)}$.
Let us fix the data on $M$:
$$\{ p_0, p_{\infty}, N_0, N_{\infty},  U\}.$$
Then these estimates hold uniformly with the same constants $\epsilon$ and $C_0$, 
among all pre-admissible  bounded Hamiltonians 
with bounds $||f||_{C^{l+1}(M)} \leq C$ by a constant $C$,
where $l$ is the fixed degree on the Sobolev  space we used.
\end{remark}

\vspace{3mm}

\begin{proof}
Let us verify 
proposition  \ref{C1-bound}.
We proceed  by contradiction argument.
Suppose contrary, and choose a sequence
$u_i \in {\frak C}_i$ with $||du_i||_{C^0({\mathbb R} \times S^1)} \to \infty$.

{\bf Step 1:}
By applying lemma $3.3$ in [HV] to $X =M_i$,  
one can find $\epsilon_i \to 0$ and $x_i \in {\mathbb R} \times S^1$ with:
 $$|du_i(x_i)| \epsilon_i \to \infty, \quad
|du_i(x)| \leq 2|du_i(x_i) | \quad (|x_i -x| \leq \epsilon_i).$$
Let us put:
$$R_i = |du_i(x_i)| \epsilon_i, \quad
v_i(x) \equiv u_i(x_i + |du_i(x_i)|^{-1} x).$$
Then $v_i$ satisfy the equation:
$$\frac{\partial}{\partial s} v_i +J(v_i) \frac{\partial}{\partial t} v_i
+ |du_i(x_i) |^{-1} \lambda_i\nabla f_i(v_i)=0$$
with $ |dv_i(0) |=1$ and $  |dv_i(x) | \leq 2$ for $x \in B_{R_i}(0) \subset {\mathbb R}^2$.
Let us choose another sequence $S_i \leq R_i$ with  $S_i \to \infty$,
so that:
 $$\text{ vol } (B_{S_i}(0) )^{\frac{1}{2}}
 |du_i(x_i) |^{-1} \lambda_i | |\nabla f_i||_{C^0(M_i)} \to 0$$
 holds.
 By elliptic regularity, 
the uniform estimate
$||dv_i||_{L^2(B_{S_i}(0))} \geq \delta$ holds for some constant $\delta >0$.

Now  we have the point-wise equalities:
\begin{align*}
& |dv_i|^2 = \omega(\frac{\partial}{\partial s}v_i, J \frac{\partial}{\partial s}v_i)
+ \omega(\frac{\partial}{\partial t} v_i, J \frac{\partial}{\partial t} v_i) \\
& =  2 \omega(\frac{\partial}{\partial s} v_i, \frac{\partial}{\partial t} v_i)
- \omega(\frac{\partial}{\partial s} v_i, J  |du_i(x_i) |^{-1} \lambda_i \nabla f_i(v_i)) \\
& \qquad \qquad
-  \omega(\frac{\partial}{\partial t} v_i, |du_i(x_i) |^{-1} \lambda_i \nabla f_i(v_i)).
\end{align*}
Then 
for sufficiently large $i >>1$, we have the estimates: 
\begin{align*}
& \int_{B_{S_i}(0)}
 |\omega(\frac{\partial}{\partial s}
 v_i,  J  |du_i(x_i) |^{-1} \lambda_i \nabla f_i(v_i))| \\
&  \leq ||dv_i||_{L^2(B_{S_i}(0)) }\text{ vol } (B_{S_i}(0) )^{\frac{1}{2}}
 |du_i(x_i) |^{-1} \lambda_i ||\nabla f_i||_{C^0(M_i)} \\
&  \leq \frac{\delta}{3} ||dv_i||_{L^2(B_{S_i}(0))} \leq  \frac{1}{3} ||dv_i||^2_{L^2(B_{S_i}(0))}.
 \end{align*}
 A similar estimate also holds:
$$ \int_{B_{S_i}(0)}
 |\omega(\frac{\partial}{\partial t}
 v_i,    |du_i(x_i) |^{-1} \lambda_i \nabla f_i(v_i))| 
  \leq  \frac{1}{3} ||dv_i||^2_{L^2(B_{S_i}(0))}.$$
Combining with these estimates, the following  must hold:
$$\delta^2 \leq  ||dv_i||^2_{L^2(B_{S_i}(0))} \leq 6 \int_{B_{S_i}(0)} v_i^*(\omega).$$

{\bf Step 2:}
Let us verify that the integral $ \int_{B_{S_i}(0)} v_i^*(\omega)$
is uniformly bounded from above.
In fact the equality holds:
$$\omega( \frac{\partial}{\partial s}v_i,  \frac{\partial}{\partial t}v_i)
 = |du_i(x_i)|^{-1} \lambda_i \frac{df_i}{ds}(v_i) 
+ |J(v_i) \frac{\partial}{\partial t}v_i +  |du_i(x_i)|^{-1} \lambda_i \nabla f_i|^2.$$
The integral of  the first term on the right hand side satisfies the estimate:
\begin{align*}
\int_{s_0}^{s_1} 
 |du_i(x_i)|^{-1} \lambda_i \frac{df_i}{ds}(v_i)ds  & 
= |du_i(x_i)|^{-1} \lambda_i (f_i(v_i(s_1, t)) - f_i(v_i(s_0, t))) \\
 & \leq 2  |du_i(x_i)|^{-1} \lambda_i || f_i||_{C^0(M_i)}.
 \end{align*}
   In particular it follows from the estimate:
    $$|\int_{B_{S_i}(0)}  |du_i(x_i)|^{-1} \lambda_i \frac{df}{ds}(v_i) |
    \leq  2 S_i  |du_i(x_i)|^{-1} \lambda_i || f_i||_{C^0(M_i)} $$
that the left hand side   is uniformly bounded from above.

Since the equality:
$$< u_i, \omega> - \lambda_i(\sup f -\inf f) = 
\int_{{\mathbb R} \times S^1}  |J(u_i) \frac{\partial}{\partial t}u_i
 +  \lambda_i \nabla f_i|^2$$
holds, the  right hand side is uniformly bounded, and hence
 the integral
$\int_{B_{S_i}(0)}  |J(v_i) \frac{\partial}{\partial t}v_i +  |du_i(x_i)|^{-1} \lambda_i \nabla f_i|^2$
is also uniformly bounded from above.

This verifies the claim.

\vspace{2mm}

\begin{remark}
One can give another argument
by   [HV] 
by use of  the basic results of holomorphic curves,
if we assume quasi-transitivity on $M$.
\end{remark}

\vspace{2mm}

{\bf Step 3:}
$L^2$ norm $||dv_i||^2_{L^2(B_{S_i}(0))} $ is  uniformly bounded from above
by step $1,2$.  
Let us choose  some $0 < \mu_i \to 0$ so that 
(1) $(\mu_i s_i, \mu_i r_i) = (1.5,1)$
for some $ 0 < r_i < s_i \leq S_i$ and
(2) $||dv_i||_{L^2(B_{s_i}(0) \backslash B_{r_i}(0) )} \to 0$.

Let us put $v_i' : B_{1.5}(0)  \to M_i$ by $v_i'(s,t) = v_i(\mu_i^{-1}s, \mu_i^{-1}t)$.
 $L^2$ norms are preserved  under  rescaling
$ ||dv_i'||_{L^2 (B_{1.5}(0) \backslash B_1(0) )}
 = ||dv_i||_{L^2(B_{s_i}(0) \backslash B_{r_i}(0) )}$, and hence
 the left hand side converges to zero.
Combining with the Cauchy-Schwartz,  
 $||dv_i'||_{L^1 (B_{1.5}(0) \backslash B_1(0) ) }$ approaches to zero. In particular
 there is some $1\leq a_i' \leq 1.5$ such that
 $$||dv_i'||_{L^1(S^1(a'_i))} \to 0$$
 where $S^1(a'_i) = \{ z \in {\mathbb C} : |z| = a'_i\}$.
 Let us put $a_i =\mu_i^{-1}a_i'$. Then the length of $v_i(S^1(a_i))$ goes to zero.

So one can cut the disk
$v_i(B_{a_i}(0))$ and 
put a small disk along $v_i(S^1(a_i))$.
By this way one  obtains a sphere $w_i $ from $v_i$, and similarly 
one  gets another sphere $w_i'$ by attaching it on  the complement
$u_i \backslash v_i(B_{a_i}(0))$ along the same boundary.

It follows from  step $1$ that 
 uniform positivity
$<w_i, \omega> \  \geq \ c >0$ must hold.
Because $<u_i, \omega> = <w_i, \omega>+<w_i', \omega> $
are the  minimal and positive number, 
$$<w_i', \omega> \  \leq \ 0$$ must be non positive. 
On the other hand there are decreasing constants $0 < \delta_i \to 0$ such that:
$$ < w'_i, \omega>  \ >  \lambda_i(\sup f -\inf f)  +
 \int_{{\mathbb R} \times S^1 \backslash B_{\epsilon_i}(x_i)}  
 | J(u_i) \frac{\partial}{\partial t}u_i +  \lambda_i \nabla f_i|^2  -  \delta_i \geq - \delta_i$$
must hold. 

{\bf Step 4:}
The above estimate implies that 
$(1)$ $< w'_i, \omega> =0$,
$(2)$ there could occur at most one bubbling, and 
$(3)$ $\int_{{\mathbb R} \times S^1 \backslash B_{\epsilon_i}(x_i)}  
 | J(u_i) \frac{\partial}{\partial t}u_i + \lambda_i \nabla f_i|$
  goes to $0$.
 
In particular 
   $s(u_i) = s_{\infty}(u_i) -s_0(u_i) \geq c >0$ must be uniformly bounded from below,
   since bubbling can occur only at one point.
$x_i \in {\mathbb R} \times S^1$ is  contained in 
one of  $[s_{\infty}(u_i), \infty) \times S^1$ or 
$(- \infty ,  s_0(u_i)] \times S^1$
or $[s_0(u_i), s_{\infty}(u_i)] \times S^1$.

Let us consider the first case.
The restriction
$u_i|(- \infty ,  s_0(u_i)] \times S^1$ must have 
uniformly bounded one derivatives.

By the construction, 
$u_i$ is holomorphic on $(- \infty, s_0(u_i)] \times S^1$ 
whose images are contained in $N_0 \subset M $. 
There is some $\delta >0$ such that $N_0$
 contains $2 \delta >0$ ball with the center $p_0$.
Then
$u_i(- \infty, S^1) =p_0$ and $d(u_i(r_i, y_i) , p_0) = \delta >0$
must hold 
for some $y_i \in S^1$ and $r_i < s_0(u_i)$.

Let us denote by $D^2(b) \subset S^2$  the disk with the radius $b$.
One may assume $ (- \infty, r_i] \times S^1 =  D^2(1)
\backslash 0  \subset S^2$ by use of translation if necessarily,
where we disregard the normalization condition
$\int_{D(1)} u^*(\omega) = \frac{1}{2}<u,\omega>$.
 Then let us choose    $  a_i >0$ with
 $ (- \infty, s_0(u_i)] \times S^1 =  D^2(1+a_i) 
\backslash 0  \subset S^2$.

We claim that $a_i>0$ is uniformly bounded from below.
In fact there is   $x_i \in S^1$ with
 $d(u(s_0(u_i), x_i), p_0) \geq 2 \delta$, and 
 the estimate $d(u_i(r_i, x_i), p_0) \leq  \delta$ must hold.
 Thus we have:
 $$d(u_i(r_i,x_i), u(s_0(u_i), x_i)) \geq \delta.$$
If $|r_i-s_0(u_i)| $ is small, then it would contradict to the assumption of
 uniform bound of $||du||$ on $(- \infty ,  s_0(u_i)] \times S^1$.
This verifies the claim.

But this would be impossible,
 since:
 $$ \int_{(-\infty, r_i] \times S^1}  | \frac{\partial}{\partial s}u_i |^2= \int_{(-\infty, r_i] \times S^1}  | \frac{\partial}{\partial t}u_i|^2 $$
   goes to $0$ by $(3)$ above, which contradicts to sublemma \ref{SU}.
   
   The second case can be considered similarly.
   
   Suppose $x_i$ is contained in $[s_0(u_i), s_{\infty}(u_i)] \times S^1$.
   Then at least one of  $[s_0(u_i) + \frac{s(u_i)}{2}, s_{\infty}(u_i)] \times S^1$
  or  $[s_0(u_i), s_0(u_i) + \frac{s(u_i)}{2}] \times S^1$
     must contain   $x_i$ for some  infinite number of $i$.
  In the former case,
one can repeat the above argument 
   over  $(- \infty, s_0(u_i)] \times S^1$.
   
   The rest case can be considered similarly.
   \end{proof}
   
\vspace{3mm}

\begin{proof}
Let us verify  lemma \ref{su-type lem}.
 By proposition  \ref{C1-bound}, 
 uniform bound:
$$||du||_{C^0({\mathbb R}  \times S^1)} \ 
 \  \leq \  C_0 < \infty$$ holds.
Let us verify 
 that
there is positive $\epsilon >0$ determined by 
$[(M_i, \omega_i, J_i)]$   with the uniform estimates:
$$ \int_{ (- \infty, s_0(u)] \times S^1} u^*(\omega) ,  \ \ 
 \int_{[s_{\infty}(u), \infty) \times S^1} u^*(\omega) \ \  \geq \epsilon.$$
 We  only verify the  estimate for the former. The latter follows by  the same argument.

{\bf Step 1:}
Let us choose $\delta >0$ so that $2\delta$ ball 
$B_{2\delta}(p_0)$ with the center $p_0$ is contained in $N_0$.
Let us put
 $s_0'(u) =  \sup \{ s \in {\mathbb R}:
\ u(( - \infty , s) \times S^1)  \subset B_{\delta}(p_0) \} $.
 $u'(s, x) \equiv u(s+ \alpha, x)$ 
still satisfies uniformity 
$||du'||_{C^0({\mathbb R}  \times S^1) }
   \leq \  C_0 < \infty$.
      By translation, 
assume $ (- \infty, s_0'(u)] \times S^1 =  D^2(1)
\backslash 0  \subset S^2$  and put
 $ (- \infty, s_0(u)] \times S^1 =  D^2(1+a) 
\backslash 0  \subset S^2$.
 $a>0$ is uniformly bounded from below
 by proposition \ref{C1-bound}.

 {\bf Step 2:}
 Notice that  $u$ is  holomorphic on $D(1+a)$, 
 and $a>0$ is uniformly bounded from below by step $1$.
Suppose $\int_{D(1+a)} u^*(\omega) < \epsilon$ could hold for 
 small $\epsilon >0$.
 Then by sublemma \ref{SU}, 
  the uniform estimates of the derivative: 
 $$ |du|  (m) \leq C \sqrt{\epsilon} $$ hold on 
 all points  $m \in D^2(1)$. 
 This is a contradiction if $\epsilon >0$ is small, since  the distance
  $d(p_0, u(s_0'(u),y)) $ attains $ \delta$ at some $y \in S^1$.
  \end{proof}
  \vspace{3mm}

\section{ Structure of the cobordism}

\subsection{Proper compactness}
Here we verify theorem \ref{product}.

\begin{sublemma}\label{convergence}
Let   $[(M_i, \omega_i, J_i)]$ be a  quasi-transitive almost K\"ahler sequence,
and suppose $  {\frak M}[(M_i,  \omega_i, J_i)]$  is  compact with respect to a minimal class.

If   a family $\{u_k \}_k \subset  {\frak D}$ satisfies
uniform bound
$||du_k||_{C^0(S^2)} \leq C$ with:
$$\lim_{k \to \infty} ||\bar{\partial}_J(u_k)|| =0$$
then a subsequence converges to some elements in 
${\frak M} [(M_i,  \omega_i, J_i)] $.
\end{sublemma}

\begin{proof}
It follows from  quasi-transitivity 
and the assumption on the uniform bound
  that for any small $\mu >0$, 
 there is
 $l_0$ and a family of 
 automorphisms $A_k$ on $M$  such that
  the estimates hold for all $k $:
$$d( \text{ im} \  A_k(u_k) , M_{l_0}) < \mu.$$

 Let $\pi : U_{\epsilon}(M_{l_0}) \to M_{l_0}$ be the holomorphic projection,
 and consider the compositions:
 $$\pi \circ A_k(u_k): S^2 \to M_{l_0}.$$
Then the Sobolev norms
$||\bar{\partial}(  \pi \circ A_k(u_k)) ||= ||(\pi \circ A_k)_*\bar{\partial} (u_k) ||
$ converge to $0$ over $M_{l_0}$.
So there is a holomorphic curve $u \in {\frak M}[(M_i, \omega_i, J_i)]$
so that  $A_k(u_k)$ are contained in a small neighborhood of 
$u$ in ${\frak B}$ for all large $k$.

 The  family $\{A^{-1}_k(u)\}_k $
 must be contained in $ {\frak M}[(M_i, \omega_i, J_i)]$,
 since $A_k$ are automorphisms which preserves $ M_0$.
 In particular $u_k$ themselves must be contained in
 a small neighborhood of $ {\frak M}[(M_i, \omega_i, J_i)]$.
 So the conclusion follows since $ {\frak M}[(M_i, \omega_i, J_i)]$  is compact.
\end{proof}

\begin{proof}
Let us verify 
 theorem \ref{product}.
For the proof we apply   analysis of
moduli theory  in section $3$.

\begin{remark}\label{lambda-est}
 The proof below verifies that 
$\lambda_0$ can be chosen uniformly 
among  pre-admissible and  bounded Hamiltonians 
with uniform bound
$||f||_{C^{l+1}(M)} \leq C$ by a constant $C$,
where $l$ is the derivative order of the Sobolev spaces
(cf. remark  \ref{unif.const}).
\end{remark}

{\bf Step 1:} We claim that there are $\lambda_0 >0$ and $C$ so that 
any $u_i \in {\frak C}_i(0 \leq \lambda \leq \lambda_0)$
satisfy uniform bounds
$||du_i||_{C^0(S^2)} \leq C$.

By proposition \ref{C1-bound}, uniformity 
$||du_i||_{C^0({\mathbb R} \times S^1)} \leq C$ holds.
It follows from the argument in the proof of lemma \ref{key lemma} that 
there is  positive $\lambda_0 >0$ such that
uniform bounds $s(u_i) =s_{\infty}(u_i) - s_0(u_i) \leq c$
hold for some 
 constant $c >0$
for all $u_i \in {\frak C}(\lambda_i)$ with $0 \leq \lambda_i \leq \lambda_0$.
In fact otherwise one would find $\lambda_i \to 0$ and $u_i \in {\frak C}(\lambda_i)$
with $s(u_i) \to \infty$. It leads us to find non trivial two spheres,
 which cannot happen.

Our claim follows, if 
 both $\{s_0(u_i)\}_i$ and $\{s_{\infty}(u_i)\}_i$ are 
uniformly bounded from both sides.

Suppose $s_0(u_i) \to - \infty$ could hold.
Then $s_{\infty}(u_i) \leq 0$ hold for all large $i$.
Let $a(s)$ be in the proof of lemma \ref{key lemma}.
Because $a(s)$ is monotone increasing,
we have the estimates:
$$a(s) \leq a(0) = \frac{1}{2}<u_i, \omega> - \lambda_i (\sup f - \inf f)$$
for all $s_{\infty}(u_i) \leq s \leq 0$,
since $u_i([s_{\infty}(u_i), \infty) \times S^1)$ is contained in $N_{\infty}$.
Notice that $u_i$ is holomorphic on the region.

There is  some $s_{\infty}(u_i) \leq s_i \leq 0$
so that the diameter of $u_i(s_i, \quad)$ goes to zero, by proposition \ref{C1-bound}.
By cutting along the circle  and putting discs on, 
one obtains two spheres $v_i^1$ and $v_i^2$ with
$p_0 \in v_i^1$ and $p_{\infty} \in v_i^2$.
Because $a(s)$ is monotone increasing, 
it follows from lemma \ref{su-type lem} that 
uniform positivity $0 < \frac{\epsilon}{2} \ \leq \ <v_i^1, \omega>$ 
must follow for all large $i$.
However this would contradict to the bounds
$<v_i^2, \omega>  \ \geq \ \frac{1}{2}<u_i, \omega> - \delta_i$
where $\delta_i \to 0$ as $i \to \infty$.

Another case can be considered similarly.

{\bf Step 2:}
We  verify that 
there is a small $\lambda_0 >0$ such that homeomorphism:
$${\frak C}_i(0 \leq \lambda <\lambda_0) \ 
\cong \ {\frak M}(M_0, \omega_0, J_0) \times [0, \lambda_0)
= {\frak M}[(M_i, \omega_i, J_i)] \times [0, \lambda_0)$$
holds for any sufficiently large $i$.

Let us
 take $u \in {\frak M}[(M_i, \omega_i, J_i)]$ and consider
 $\hat{U}(u)$ as in \ref{Hilbert-compl}.
Let $f: M \to [0, \infty)$ be a pre-admissible bounded Hamiltonian. Then 
we have the corresponding  functional 
$\hat{\frak F} : \hat{U}(u) \to \hat{\frak E}$ as in \ref{f.d.case}.
Let us put the cobordism:
\begin{align*}
{\frak C}= 
\{  \  (u', \lambda)  \in \hat{U} (u) \times [0, \infty) :  \ \
& \bar{\partial}(u') + \lambda \hat{\frak F}(u') =0 , \\
& \lambda  \geq 0,  \ \  u \in {\frak M}[(M_i, \omega_i, J_i)]  \ \}.
\end{align*}
It follows from  theorem \ref{basic-thm}
 that $ {\frak M}[(M_i, \omega_i, J_i)]$
coincides with $ {\frak M}(M_0, \omega_0, J_0)$, and 
is strongly regular and compact. 
So
the embedding: 
$${\frak C}(0 \leq \lambda <\lambda_0) \ 
\supset \ {\frak M}(M_0, \omega_0, J_0) \times [0, \lambda_0)$$
holds for some positive $\lambda_0 >0$,
by the infinite dimensional implicit function theorem.

One can make a small perturbation of $f$ 
with $df|M_i = d(f|M_i)$ for all sufficiently large $i$
by properness, so that 
${\frak C}(0 \leq \lambda <\lambda_0)$ contain 
${\frak C}_i(0 \leq \lambda <\lambda_0)$ regularly.
So we obtain the embeddings:
$${\frak C}_i(0 \leq \lambda <\lambda_0) \ 
\supset \ {\frak M}(M_0, \omega_0, J_0) \times [0, \lambda_0).$$

If for any small $\lambda_0 >0$, there could exist some $i$ such that
${\frak C}_i(0 \leq \lambda <\lambda_0) \ 
= \ {\frak M}(M_0, \omega_0, J_0) \times [0, \lambda_0)$ 
were not satisfied,
 then  there should exist divergent sequence
$u_i \in   {\frak C}_i(0 \leq \lambda <\delta_i)$ with $\delta_i \to 0$.
In particular $\lim_i ||\bar{\partial}_J(u_i)||=0$ must hold, 
which cannot happen by sublemma \ref{convergence} with step $1$.
\end{proof}

\subsection{A new inequality}
Let us induce a new inequality which arises from comparison between structure of 
cobordism of the moduli space and of iteration of  Hamiltonian diffeomorphisms.
Let $[(M_i, \omega_i, J_i)]$ be an almost K\"ahler sequence, and fix the data 
$\{ p_0, p_{\infty}, N_0, N_{\infty},  U\}$ on $M$.

 \vspace{3mm}
 
 \begin{proof}
 Let us verify 
  theorem \ref{inequ}.
  
 {\bf Step 1:}
 We verify that there is a constant $C$ so that 
 for any $f$ with $||f||_{C^{l+1}(M)} \leq 1$,
 the inequality $C \leq \text{ Cob }(f) \text{ As}(f) $ holds, where
 $l$ is the order on the Sobolev derivatives (see \ref{f.d.pre}).
 
 Recall a uniform property in remark \ref{unif.const}.
 It follows from the proof of lemma \ref{key lemma} that 
  there is $\epsilon >0$ so that  any $f$ as above satisfies  the uniform estimate
 $|\delta(\tilde{F}_i , m_i)| > \epsilon$,
 where $\tilde{F}_i$ are the Hamiltonian diffeomorphisms with respect to $\lambda_i f_i$.

It follows from  \ref{proof of unif.cyc} in the proof of theorem \ref{unif.cyc} step $4$
that there is a constant $C$ such that  the estimate:
$$ n \epsilon \leq C\lambda  ||d(F_i)^n||_{C^0(M_{i})}$$
holds for any  $n \geq 1$,
where $\lambda = \lim_{i \to \infty} \lambda_i$.
So the lower bound holds:
  $$ \frac{\epsilon}{C \lambda} 
  \leq \text{ As} (f)  . $$

  {\bf Step 2:}
 Theorem \ref{product} and remark \ref{lambda-est}
 verifies uniform positivity:
 $$\text{ Cob }(f) \geq \lambda_0 >0.$$
  $\lambda_i$ is chosen so that the estimate
 $\lambda_i  > \lambda_0$ holds in lemma \ref{top.bound}, since
 we have chosen $(u_i, \lambda_i)  \in {\frak C}_i(\lambda_i)$
 so that the family $\{u_i\}_i$ diverges.

We claim that one can choose $\lambda_i$
 so that $\lambda_i \to \text{ Cob }(f)$ holds.
 Suppose contrary.
 By definition of
  $\text{ Cob }(f)$, the estimate  $ \lambda_i \geq 
 \text{ Cob }(f)$ holds, and so 
  $ \lambda_i >
 \text{ Cob }(f)$ has to hold. Then 
 for any $ \lambda_i > \mu > \mu' > \text{ Cob }(f)$,
 there should exist infinite number of  $i$ such that 
 ${\frak C}_i(\delta)$ are empty  or non compact for all $\mu' \leq \delta \leq \mu$.
Because 
 ${\frak C}_i$ cannot give $S^1$ freely zero cobordism to the moduli space
 of holomorphic curves in lemma \ref{top.bound} (proposition $2.7$ in [HV]), 
${\frak C}_i(0 \leq \delta \leq \mu)$ is also non compact.
Then one can choose $0 < \lambda_i \leq \mu$, which gives a contradiction.
This verifies the claim.
 
 Combining with these, we obtain 
  the uniform lower bound: 
$$\epsilon C^{-1} \leq \text{ Cob }(f) \text{ As} (f).$$

{\bf Step 3:}
Let us choose some $1 \geq \alpha >0$ with $||\alpha f ||_{C^{l+1}(M)} \leq 1$.
Clearly we have the equality:
$$\text{ Cob }(\alpha f) = \alpha^{-1} \text{ Cob }(f).$$

Let  $\tilde{F}$ be the Hamiltonian diffeomorphism with respect to $\alpha f$.
Then we have the equality
$\tilde{F}_t =F_{\alpha t}$.
Let us put $\alpha^{-1} m = k_m - \alpha_m$ with $k_m \in {\mathbb N}$ and $0 \leq  \alpha_m < 1$.
Then:
\begin{align*}
 \frac{1}{m}   ||d (F_{i})^{ m}||_{C^0(M_{i})} 
 & =  \frac{1}{m} ||d(\tilde{F}_{i} )_{\alpha^{-1} m}||_{C^0(M_{i})}   \\
& = \frac{1}{m} ||d(\tilde{F}_{i})_{- \alpha_m} \circ d(\tilde{F}_{i} )^{k_m}||_{C^0(M_{i})} \\
 & \geq C_0    \frac{\alpha^{-1}}{ \alpha^{-1} m} ||d(\tilde{F}_{i} )^{k_m}||_{C^0(M_{i}) } \\
 & \geq  C_0  \frac{\alpha^{-1}}{k_m} ||d(\tilde{F}_{i} )^{k_m}||_{C^0(M_{i})}.
\end{align*}
 By letting $m,i  \to \infty$, 
  we obtain the estimate:
$$\text{ As} (f) \geq C_0  \alpha^{-1} \text{ As}(\alpha f).$$

So we obtain the inequality:
\begin{align*}
\text{ Cob }(f)  \text{ As} (f)  
&  \geq \alpha \text{ Cob }( \alpha f)  C_0  \alpha^{-1} \text{ As} (\alpha f) \\
&  = C_0 \text{ Cob }( \alpha f)  \text{ As}(\alpha f) \ \geq \ C_0'.
\end{align*}
 \end{proof}

\vspace{3mm}

For any  bounded Hamiltonian $f$,  let us  introduce another  invariant:
 $$
 \widehat{  Cob}(f)  = 
 \sup_{\lambda \geq 0} \{ \ \lambda : {\frak C} (\mu) 
 \text{ are  non empty and compact  for all } 
 0 \leq \mu \leq  \lambda  \}.$$ 
 It would be interesting to compare 
 $\widehat{ Cob}(f)$  with ${ Cob}(f)$ 
 on some reasonable class of spaces such as the one in theorem \ref{product}.

\large

\vspace{1cm}

Tsuyoshi Kato

Department of Mathematics

Faculty of Science

Kyoto University

Kyoto 606-8502

Japan

\vspace{5mm}

\end{document}